\input ./style/arxiv-general.cfg
\documentclass[aap,MSNbibl,seceqn,dvips]{arximspdf}
\makeatletter
   \@ifpackageloaded{graphicx}{}{\usepackage{graphicx}}
\makeatother

\doi{10.1214/15-AAP1119}
\volume{26}
\issue{3}
\pubyear{2016}
\firstpage{1362}
\lastpage{1406}
\docsubty{FLA}

\makeatletter
\newcommand{\fontas}{\fontsize{8.36}{10.36}\selectfont}
\newcommand{\rrVert}{\Vert}
\newcommand{\llVert}{\Vert}
\renewcommand{\mid}{|}
\newproclaim{definition}{Definition}[section]
\newtheorem{theorem}[definition]{Theorem}
\newtheorem{lemma}[definition]{Lemma}
\newtheorem{proposition}[definition]{Proposition}
\newtheorem{corollary}[definition]{Corollary}
\newproclaim{fact}[definition]{Fact}
\newcommand\G{\mathbf{H}}
\newcommand{\Zkc}{Z_{\beta}}
\newcommand{\Zkb}{Z_{\beta,\mathrm{bal}}}
\newcommand{\Zkg}{Z_{\beta,\mathrm{tame}}}
\newcommand{\pr}{\mathbb P}
\newcommand\core{\mathrm{core}}
\newcommand\back{\mathrm{back}}
\newcommand\rest{\mathrm{rest}}
\newcommand\free{\mathrm{free}}
\newcommand\Forb{\operatorname{Forb}}
\newcommand\Bal{\operatorname{Bal}}
\newcommand\cA{\mathcal{A}}
\newcommand\cB{\mathcal{B}}
\newcommand\cH{\mathcal{H}}
\newcommand\eul{\mathrm{e}}
\newcommand\eps{\varepsilon}
\newcommand\Erw{\mathbb{E}}
\newcommand{\vecone}{\mathbf{1}}
\newcommand{\Vol}{\operatorname{Vol}}
\newcommand{\Bin}{\operatorname{Bin}}
\newcommand{\bink}[2]{{{#1}\choose{#2}}}
\newcommand{\binkd}[2]{\pmatrix{#1\cr #2}}

\newcommand{\brk}[1]{[{#1}]}
\newcommand{\brkk}[1]{\bigl[{#1}\bigr]}
\newcommand{\brkkk}[1]{\biggl[{#1}\biggr]}
\newcommand{\brkkkkp}[1]{\Biggl[{#1}\Biggr]}

\newcommand{\scal}[2]{\langle{#1},{#2} \rangle}
\newcommand\RR{\mathbb{R}}
\makeatother

\begin{document}
\begin{frontmatter}

\title{A positive temperature phase transition in random hypergraph
2-coloring\thanksref{T1}}
\runtitle{A positive temperature phase transition}

\begin{aug}
\author[A]{\fnms{Victor}~\snm{Bapst}\ead[label=e1]{bapst@math.uni-frankfurt.de}},
\author[A]{\fnms{Amin}~\snm{Coja-Oghlan}\corref{}\ead[label=e2]{acoghlan@math.uni-frankfurt.de}}
\and
\author[A]{\fnms{Felicia}~\snm{Ra\ss mann}\ead[label=e3]{rassmann@math.uni-frankfurt.de}}
\runauthor{V. Bapst, A. Coja-Oghlan and F. Ra\ss mann}
\affiliation{Goethe University}
\address[A]{Mathematics Institute\\
Goethe University\\
10 Robert Mayer St\\
Frankfurt 60325\\
Germany\\
\printead{e1}\\
\phantom{E-mail: }\printead*{e2}\\
\phantom{E-mail: }\printead*{e3}}
\end{aug}
\thankstext{T1}{Supported by the European Research Council under the
European Union's Seventh Framework Programme (FP/2007-2013)/ERC Grant
Agreement no. 278857--PTCC.}

%
\received{\smonth{10} \syear{2014}}
%
\revised{\smonth{3} \syear{2015}}

%
\begin{abstract}
Diluted mean-field models are graphical models in which the geometry of
interactions
is determined by a sparse random graph or hypergraph. Based on a
nonrigorous but
analytic approach called the ``cavity method'', physicists have
predicted that in
many diluted mean-field models a phase transition occurs as the inverse
temperature
grows from $0$ to $\infty$
[\textit{Proc. National Academy of Sciences} \textbf{104} (2007) 10318--10323].
In this paper, we establish
the existence
and asymptotic location of this so-called condensation phase transition
in the random hypergraph $2$-coloring problem.
\end{abstract}

%
\end{frontmatter}

\section{Introduction and results}

\subsection{Background and motivation}

Statistical mechanics models of ``disordered system'' such as glasses
or spin-glasses are notoriously difficult to study analytically.
Nonetheless, since the early 2000s physicists have developed an
analytic but nonrigorous approach, the so-called \textit{cavity
method}, to
put forward precise conjectures on an important class of models called
\textit{diluted mean-field models}.
These are models where the geometry of interactions between individual
``sites'' is determined by a sparse random graph or hypergraph.
Apart from models of inherent physical interest, the cavity method has
since been applied
to a wide variety of problems in combinatorics, computer science,
information theory and compressive sensing~\cite{MM,KMSSZ}.
What these problems have in common is that there are ``variables'' and
``constraints'' whose mutual interaction is governed by a sparse random
hypergraph.
In effect, it has become an important research endeavour to provide a
rigorous mathematical foundation for the cavity method.
The present paper contributes to this effort.

Among the various predictions deriving from the cavity method, perhaps
the most intriguing ones pertain to the existence and location of phase
transitions.
In particular, according to the cavity method in a variety of models
there occurs a so-called \textit{condensation phase transition}.
This is a phenomenon that is ubiquitous in physics.
Its role in the context of structural glasses goes back to the work of
Kauzmann in the 1940s~\cite{Kauzmann48}.
However, there are but a few rigorous results on the condensation phase
transition in diluted mean-field models.

The aim of the present work is to establish the existence and
asymptotic location of the condensation phase transition
in a well-studied diluted mean-field model, the \textit{random hypergraph
$2$-coloring problem}.
To define this model, we recall that a $k$-uniform hypergraph $H$
consists of a finite set $V_H$ of vertices and a set $E_H$ of edges,
which are subsets of $V_H$ of size $k$.
For a $k$-uniform hypergraph $H$ and a map $\sigma:V_H\rightarrow
\{{-1,1} \}$
we let $E_H(\sigma)$ be the number of edges
$e\in E_H$ such that $\vert\sigma(e)\vert=1$, that is, either all
vertices of
$e$ are set to $1$ or to $-1$ under $\sigma$.
Thus, if we think of $\sigma$ as a coloring of the vertices of $H$
with two colors, then $E_H(\sigma)$ is the number of monochromatic edges.
The Hamiltonian $E_H$ gives rise to a Boltzmann distribution $\pi
_{H,\beta}$ on the set of all maps $\sigma:V_H\rightarrow\{
{-1,1} \}$ in the
usual way: we let
%
%
\begin{eqnarray}\label{eqBoltzmann}
\pi_{H,\beta}[\sigma]=\frac{\exp(-\beta
E_H(\sigma))}{Z_\beta
(H)}
\nonumber\\[-8pt]\\[-8pt]
\eqntext{\displaystyle\mbox{where } Z_\beta(H)=\sum_{\tau:V_H\rightarrow
\{{-1,1} \}}\exp\bigl(-\beta
E_H(\tau)\bigr)}
\end{eqnarray}
is the partition function.
We refer to $\beta$ as the \textit{inverse temperature}.
Clearly, as \mbox{$\beta\rightarrow\infty$} the Boltzmann distribution $\pi
_{H,\beta}$ will place more and more weight on maps $\sigma$ with
fewer and fewer monochromatic edges.
For a given hypergraph $H$, the key object of interest is the function
$\beta\mapsto\frac{1}n\ln Z_\beta(H)$, the \textit{free entropy}.

While the definition~(\ref{eqBoltzmann}) makes sense for any
hypergraph $H$,
in the diluted mean-field model
the hypergraph itself is random.
More specifically, we consider the random hypergraph $H_k(n,p)$ on $n$
vertices $V= \{{1,\ldots,n} \}$,
in which each of the $\bink nk$ possible hyperedges comprising of $k$
distinct vertices is present with probability $p\in[0,1]$ independently.
Throughout the paper, we always let $\beta\in[0,\infty)$ and
$p=d/\bink{n-1}{k-1}$, where $d>0$ is a real number and $k\geq3$ is
an integer.
The parameters $d,k$ and $\beta$ are going to remain fixed while we
are going to let $n\rightarrow\infty$.
The main objective is to determine
%
\begin{equation}
\label{eqFreeEntropyDensity} \Phi_{d,k}(\beta)=\lim_{n\rightarrow
\infty}
\frac{1}n\Erw\brkk{\ln Z_\beta\bigl(H_k(n,p)\bigr)},
\end{equation}
the \textit{free entropy density}.
Of course, in~(\ref{eqFreeEntropyDensity}) the expectation is over the
choice of the random hypergraph $H_k(n,p)$.

An obvious question is whether the limit~(\ref{eqFreeEntropyDensity})
exists for all $d,k,\beta$.
That this is indeed the case follows from an application of the
combinatorial interpolation method from~\cite{bayati}.
Furthermore, a standard application of Azuma's inequality shows that
for any $d,k,\beta$ the sequence
$\{\frac{1}n\ln Z_\beta(H_k(n,p))\}_{n}$ converges to $\Phi
_{d,k}(\beta
)$ in probability.

\subsection{The main result}

In this paper, we establish the existence and approximate location of
the condensation phase transition in random hypergraph $2$-coloring.
More specifically, we are going to obtain a formula that determines the
location of the condensation phase transition
up to an error $\eps_k$ that tends to $0$ for large $k$.
This is the first (rigorous) result that determines the condensation
phase transition within such accuracy in terms of the finite parameter
$\beta$ (the ``positive temperature'' case, in physics jargon).

We call $\beta_0>0$ \textit{smooth} if there exists $\eps>0$ such that
the function
$\beta\in(\beta_0-\eps,\beta_0+\eps)\mapsto\Phi_{d,k}(\beta)$
admits an expansion as an absolutely convergent power series around
$\beta_0$.
Otherwise, we say that a \textit{phase transition} occurs at $\beta_0$.
With these conventions, we have the following theorem.

%
\begin{theorem}\label{Thmcond}
For any fixed number $C>0$, there exists a sequence $\eps_k>0$ with
$\lim_{k\rightarrow\infty}\eps_k=0$ such that
the following is true.
Let
\[
\Sigma_{k,d}(\beta) = ( \beta+1 ) \exp(- \beta+ k \ln2 )\ln2 -2
\biggl( \frac{d}{k}- 2^{k-1} \ln2 + \ln2 \biggr).
\]
\begin{longlist}[(iii)]
\item[(i)] If $d/k < 2^{k-1} \ln2 - \ln2 - \eps_k$, then any $\beta>0$
is smooth and
%
%
\begin{equation}
\label{eqThmcondFreeEntropy} \Phi_{d,k}(\beta)=\ln2+\frac{d}k\ln
\bigl(1-2^{1-k}\bigl(1-\exp(-\beta)\bigr)\bigr).
\end{equation}
\item[(ii)] If $2^{k-1} \ln2 - \ln2+ \eps_k <d/k<2^{k-1} \ln2+C$, then
$\Sigma_{k,d}(\beta)$ has a unique zero $\beta_c(d,k) \geq k \ln2 $ and:
\begin{itemize}
\item any $\beta\in(0, \beta_c(d,k) + \eps_k)$ is smooth and $\Phi
_{d,k}(\beta)$ is given by~(\ref{eqThmcondFreeEntropy}),
\item there occurs a phase transition at $\beta_c(d,k) + \eps_k$,
\item for $\beta>\beta_c(d,k) + \eps_k$ we have
\[
\Phi_{d,k}(\beta)<\ln2+\frac{d}k\ln\bigl(1-2^{1-k}
\bigl(1-\exp(-\beta)\bigr)\bigr).
\]
\end{itemize}
\end{longlist}
\end{theorem}

In summary, Theorem~\ref{Thmcond} shows that in the case that the
``density'' $d/k$ of the random hypergraph is less than about
$2^{k-1}\ln2-\ln2$, there does not occur a phase transition for any
finite $\beta$.
By contrast, for slightly larger densities there is a phase transition.
Its approximate location is given by $\beta_c(d,k)$.
While in Theorem~\ref{Thmcond} this value is determined implicitly as
the zero of $\Sigma_{k,d}(\beta)$,
it is not difficult to obtain the expansion
%
%
\begin{equation}
\beta_c(d,k)=(k-1)\ln2+\ln k + 2 \ln\ln2 - \ln c +
\delta_k,
\end{equation}
where $c = d/k - 2^{k-1} \ln2 + \ln2$ and $\lim_{k\rightarrow\infty
}\delta
_k=0$. Furthermore, the proof of Theorem~\ref{Thmcond} shows that
there exists $c_1>0$ such that
$\eps_k\leq k^{c_1}2^{-k}$.
Thus, Theorem~\ref{Thmcond} determines the critical density from that
on a phase transition starts to occur
and the critical $\beta_c(d,k)$ up to an error term that decays
exponentially with~$k$.

\subsection{Discussion and related work}
In this section, we explain how Theorem~\ref{Thmcond} relates to the
predictions based on the physicists' ``cavity method''.
We also comment on further related work.
As usual, we say that an event occurs \textit{asymptotically almost
surely} (a.a.s.) if its probability converges to $1$ as $n\rightarrow
\infty$.

\subsubsection{The ``entropy crisis''}
Theorem~\ref{Thmcond} is perfectly in line with the picture sketched by
the (nonrigorous) cavity method,
and its proof is inspired by the physicists' notion that the
condensation phase transition results from an ``entropy crisis''~\cite
{pnas,MM}.
More specifically, it is expected that already for densities much
smaller than the one treated in Theorem~\ref{Thmcond}, namely for $d/k$
beyond about $2^{k-1}\ln k/k$ and for large enough $\beta$,
the Boltzmann distribution can be approximated by a convex combination
of probability measures corresponding to ``clusters'' of $2$-colorings
a.a.s. That is, there exist sets ${{\mathcal C}}_{\beta,1},\ldots
,{{\mathcal
C}}_{\beta,N}\subset\{{-1,1} \}^n$ and small numbers $0<\eps<\delta
$ such that:
\begin{itemize}
\item if $\sigma,\tau\in{{\mathcal C}}_{\beta,i}$ for some $i$,
then $\scal
{\sigma}\tau>(1-\eps)n$,
\item if $\sigma\in{{\mathcal C}}_{\beta,i},\tau\in{{\mathcal
C}}_{\beta,j}$ with
$i\neq j$, then $\scal{\sigma}\tau<(1-\delta)n$.
\end{itemize}
Moreover, with $Z_{\beta,i}=\sum_{\tau\in{{\mathcal C}}_{\beta
,i}}\exp
(-\beta E_{H_k(n,p)}(\tau))$ the volume of ${{\mathcal C}}_{\beta,i}$,
we have
\[
\Biggl\llVert\pi_{H_k(n,p),\beta}[ \cdot]-\sum_{i=1}^N
\frac{Z_{\beta,i}}{Z_\beta(H_k(n,p))}\cdot\pi_{H_k(n,p),\beta}[
\cdot\mid{{\mathcal
C}}_{\beta,i}]\Biggr\rrVert_{\mathrm{TV}}<\exp\bigl(-\Omega
(n)\bigr).
\]
Given a hypergraph, the definition of the ``clusters'' ${{\mathcal
C}}_{\beta,i}$ is somewhat canonical (under certain assumptions); we will
formalise the construction in Section~\ref{secoutline}.

With the cluster decomposition in place, the physics story of how the
condensation phase transition comes about goes as follows.
If $\beta$ is sufficiently small,
we have $\max_{i\leq N}\ln Z_{\beta,i}\leq\ln Z_{\beta}(H_k(n,p)
)-\Omega(n)$ a.a.s. That is, even the largest cluster only captures an
exponentially small
fraction of the overall mass $Z_{\beta}(H_k(n,p))$.
Now, as we increase $\beta$ (while $d/k$ remains fixed), both
$Z_{\beta}(H_k(n,p))$ and $\max_{i\leq N}Z_{\beta,i}$ decrease.
But $Z_{\beta}(H_k(n,p))$ drops at a faster rate.
In fact, for large enough densities $d/k$ there might be a critical
value $\beta_*$ where the gap between
$\max_{i\leq N}\ln Z_{\beta,i}$ and $\ln Z_{\beta}(H_k(n,p))$ vanishes.
This $\beta_*$ should mark a phase transition.
This is because $\max_{i\leq N}\ln Z_{\beta,i}$ and $\ln Z_{\beta
}(H_k(n,p))$ cannot both extend analytically to $\beta>\beta_*$,
as otherwise we would arrive at the absurd conclusion that $\max
_{i\leq N}Z_{\beta,i}> Z_{\beta}(H_k(n,p))$.

The proof of Theorem~\ref{Thmcond} is based on turning this ``entropy
crisis'' scenario into a rigorous argument.
To this end, we establish a rigorous version of the above ``cluster
decomposition'' and, crucially, an estimate of the cluster volumes
$Z_{\beta,i}$.
The arguments that we develop for these problems partly build upon
prior work from~\cite{Barriers,AchlioptasMoore,Lenka}.

The key difference between~\cite{Barriers,AchlioptasMoore,Lenka} and
the present work is the presence of the parameter $\beta$.
More precisely, \cite{Barriers,AchlioptasMoore,Lenka} dealt with \textit
{proper} hypergraph $2$-colorings, that is, maps $\sigma:V\rightarrow
\{{-1,1} \}$ such that $E_H(\sigma)=0$.
Thus, the Boltzmann distribution in those papers is just the uniform
distribution over proper $2$-colorings, and
the partition function is the number of proper $2$-colorings.
In a sense, this corresponds to setting $\beta=\infty$ in the present setup.
In particular, the only parameter present in \cite
{Barriers,AchlioptasMoore,Lenka} is the average degree $d$ of the
random hypergraph, whereas
in the present paper we deal with a two-dimensional phase diagram
governed by $d$ and, additionally, $\beta$.
Of course, from a ``classical'' statistical physics viewpoint it seems
less natural to vary the parameter $d$ that governs the
geometry of the system and fix $\beta$ than to fix $d$ and vary $\beta$.
Theorem~\ref{Thmcond} encompasses the latter case.

To prove Theorem~\ref{Thmcond}, we extend some of the arguments
from~\cite{Barriers,AchlioptasMoore,Lenka}.
In particular, we provide a ``finite-$\beta$'' version of the second
moment arguments from~\cite{Barriers,Lenka}.
Independently of the present work, a similar extension was obtained by
Achlioptas and Theodoropoulos~\cite{AchlioptasTheodoropoulos}.
In addition, we extend the argument for estimating the cluster size
from~\cite{Lenka} to the case of finite $\beta$.
Moreover, the argument that we develop for inferring the condensation
transition from the second moment
method and the estimate of the cluster size draws upon ideas developed
for the $\beta=\infty$ case in~\cite{Barriers,bapst,Lenka}.
Especially with respect to the estimate of the cluster size, dealing
with finite $\beta$ requires substantial additional work and ideas.

\subsubsection{Prior work on condensation}

The first rigorous result on a genuine condensation phase transition in
a diluted mean field model is
due to Coja-Oghlan and Zdeborov\'a~\cite{Lenka}, who dealt with the
proper hypergraph $2$-colorings (i.e., the $\beta=\infty$ case of the
problem considered here).
Thus, the only parameter in~\cite{Lenka} is~$d$.
The main result of~\cite{Lenka} is that there occurs a condensation
phase transition
at $d/k=2^{k-1}\ln2-\ln2+\gamma_k$, where $\lim_{k\rightarrow
\infty
}\gamma_k=0$.
Up to the error term $\gamma_k$, the result confirms a prediction from
\cite{DRZ08}.
Moreover, as Theorem~\ref{Thmcond} shows, the result from~\cite{Lenka}
matches the smallest density for which a condensation phase transition
occurs for a finite $\beta$.
In this sense, \cite{Lenka} determines the intersection of the
``condensation line'' in the two-dimensional phase diagram of
Theorem~\ref{Thmcond} with the $d$-axis.
Additionally, Bapst, Coja-Oghlan, Hetterich, Ra\ss mann and
Vilenchik~\cite{bapst} determined the condensation phase transition in
the random graph coloring problem.
This is the zero-temperature case of the Potts antiferromagnet on the
Erd\H{o}s--R\'enyi random graph.
Thus, also in~\cite{bapst} the parameter $\beta$ is absent.

The only prior (rigorous) paper that explicitly deals with the positive
temperature case is the recent work of Contucci, Dommers, Giardina and
Starr~\cite{CDGS}.
They study the $k$-spin Potts antiferromagnet on the Erd\H{o}s--R\'
enyi random graph with finite~$\beta$
and show that for certain values of the average degree a condensation
phase transition exists.
But to the extent that the results are comparable, \cite{CDGS} is less
precise than Theorem~\ref{Thmcond}.
Indeed, a direct application of the approach from~\cite{CDGS} to the
present problem would determine
$\beta_c(d,k)$ only up to an additive error of $\ln k$, rather than an
error that diminishes with $k$.
This is due to two technical differences between the present work
and~\cite{CDGS}.
First, the second moment argument required in the case of the $k$-spin
Potts antiferromagnet is technically \textit{far} more challenging than
in the present case.
In effect, an enhanced version of the second moment argument along the
lines of~\cite{Lenka} (with explicit conditioning on the cluster size)
is not available in the Potts model.
Second, \cite{CDGS} employs a conceptually less precise estimate of
the cluster size than the one we derive.
More precisely, \cite{CDGS} essentially neglects the entropic
contribution to the cluster size, leading to under-estimate the typical
cluster size significantly.

The condensation line at finite $\beta$ in the Potts antiferromagnet
on the Erd\H{o}s--R\'enyi random graph was studied by
Krzakala and Zdeborov\'a~\cite{LenkaFlorent2} by means of nonrigorous
techniques.
They predict the location of the condensation line in terms of an
intricate fixed-point problem.
(While conjectured to yield the exact location of the phase transition
for large enough average degrees $d$,
no explicit expansion for large $d$ such as the one of Theorem~\ref
{Thmcond} was given.)

\section{Preliminaries and notation}

Because we take the limit $n\rightarrow\infty$ and due to the
presence of the
sequences $\eps_k,\eps_k'$, Theorem~\ref{Thmcond} is an asymptotic
statement in both $n$ and $k$.
Therefore, throughout the paper we tacitly assume that both $n,k$ are
sufficiently large.

We use the standard $O$-notation when referring to the limit
$n\rightarrow
\infty$.
Thus, $f(n)=O(g(n))$ means that there exist $C>0$, $n_0>0$ such that
for all $n>n_0$ we have $\vert f(n)\vert\leq C\cdot\vert
g(n)\vert$.
In addition, we use the standard symbols $o(\cdot),\Omega(\cdot
),\Theta(\cdot)$.
In particular, $o(1)$ stands for a term that tends to $0$ as
$n\rightarrow
\infty$. We adopt the common notation that for the symbol $\Omega
(\cdot)$ the sign matters, that is, $f(n)=\Omega(g(n))$ means that
there exist $C>0$, $n_0>0$ such that for all $n>n_0$ we have $f(n)\ge
C\cdot g(n)$ whereas $f(n)=-\Omega(g(n))$ implies $-f(n)\ge C\cdot
g(n)$ for all $n>n_0$.

Additionally, we use asymptotic notation with respect to $k$.
To make this explicit, we insert $k$ as an index.
Thus, $f(k)=O_k(g(k))$ means that there exist $C>0$, $k_0>0$ such that
for all $k>k_0$ we have $\vert f(k)\vert\leq C\cdot\vert
g(k)\vert$.
Further,\vspace*{1pt} we write $f(k)=\tilde O_k(g(k))$ to indicate that there exist
$C>0$, $k_0>0$ such that
for all $k>k_0$ we have $\vert f(k)\vert\leq k^C\cdot\vert
g(k)\vert$.
An analogous convention applies to $o_k(\cdot),\Omega_k(\cdot)$ and
$\Theta_k(\cdot)$. Notice that here as well we have $\Omega_k(\cdot
) \ne-\Omega_k(\cdot)$.\vspace*{2pt}

Throughout\vspace*{1pt} the paper, we set $p=d/\bink{n-1}{k-1}$. The \textit{degree}
of a vertex $v\in V$ in a hypergraph $H=(V,E)$ is the number of all
edges $e\in E$ that contain $v$. We let $e(H)$ denote the total number
of edges of the hypergraph~$H$.

If $L$ is an integer, then we write $[L]$ for the set $ \{
{1,\ldots,L} \}$.
Moreover, $\cH(z)=-z\ln z-(1-z)\ln(1-z)$
denotes the entropy function.
Further, we need the following instalment of the Chernoff bound.

%
\begin{lemma}[(\cite{janson}, page~29)]\label{lemchernoffsum}
Assume that $X_1,\ldots,X_n$ are independent random variables such
that $X_i$ has a Bernoulli distribution with mean $p_i$.
Let $\lambda=\Erw[X]$ and set
$\phi(x)=(1+x)\ln(1+x)-x$.
Then
\[
\pr\brk{X\ge\lambda+t}\le\exp\bigl(-\lambda\phi(t/\lambda
)\bigr), \qquad\pr
\brk{X\le\lambda-t}\le\exp\bigl(-\lambda\phi(-t/\lambda)\bigr)
\]
for any $t>0$. In particular, $\pr\brk{X\ge t\lambda} \le\exp
(-t\lambda\ln(t/\eul))$ for any $t>1$.
\end{lemma}

It is well known that $\ln Z_\beta$, the key quantity that we are
interested in, enjoys the following ``Lipschitz property''.

%
\begin{fact}\label{eqLip}
Let $H$ be a hypergraph and obtain another hypergraph $H'$ from $H$ by
either adding or removing a single edge.
Then $\vert\ln Z_\beta(H)-\ln Z_\beta(H')\vert\leq\beta$.
\end{fact}

This Lipschitz property implies the following concentration bound
for $\ln Z_\beta(H_k(n,p))$.

%
\begin{lemma}\label{LemmaZAzuma}
For any $\alpha>0$ there is $\delta=\delta(\alpha)>0$ such that
\[
\pr\brkk{\bigl\vert\ln Z_\beta\bigl(H_k(n,p)\bigr)-\Erw
\bigl[\ln Z_\beta\bigl(H_k(n,p)\bigr)\bigr]\bigr\vert>
\alpha n}<\exp(-\delta n).
\]
\end{lemma}
\begin{pf}
This is immediate from Fact~\ref{eqLip} and McDiarmid's inequality
\cite{McDiarmid}, Theorem~3.8.
\end{pf}

Throughout the paper, it will be convenient to work with two other
random hypergraph models.
More precisely, for integers $n,m>0$ we let $H_k(n,m)$ denote the random
hypergraph on the vertex set $[n]$
obtained by choosing exactly $m$ edges without replacement uniformly at
random from all possible edges, each comprising of $k$ distinct
vertices from $[n]$. This random hypergraph model will be used
essentially in Section~\ref{SecplantedCluster}. The disadvantage of
this model is the fact that the edges are not mutually independent.
Therefore, to simplify calculations in Section~\ref{Secbcrit} we let
$H'_k(n,m)$ denote the random hypergraph on the vertex set $[n]$ obtained
by choosing $m$ edges uniformly and independently at random. In this
model, we may choose the same edge more than once, however, the
following statement shows that this is quite unlikely.

%
\begin{fact}\label{LemmadoubleEdge}
Assume that $m=m(n)$ is a sequence such that $m=O(n)$ and let $\cA$ be
the event that $H'_k(n,m)$ has no multiple edges. Then $\pr\brk{\neg
\cA
}=O(1/n^{k-2})$.
\end{fact}

We relate the expected values of the partition functions of $H_k(n,m)$ and
$H'_k(n,m)$ in Section~\ref{Secfirst}.

\section{Outline}\label{secoutline}

\textit{Throughout this section let $0\leq d/k\leq2^{k-1}\ln2+O_k(1)$.}

The proof of Theorem~\ref{Thmcond} is based on establishing the
physicists' notion of an ``entropy crisis'' rigorously.
To this end, we are going to trace two key quantities.
First, the free entropy density $\Phi_{d,k}(\beta)$, which mirrors
the typical value of the partition function $Z_\beta(H_k(n,p))$.
Second, the size of the ``cluster'' of a typical $\sigma$ chosen from
the Boltzmann distribution.
More specifically, we are going to argue that it is sufficient to study
the (appropriately defined) ``cluster size'' in a
certain auxiliary probability space, the so-called ``planted model''.
Ultimately, it will emerge that the condensation phase transition marks
the point where the cluster size in the
planted model equals the typical value of $Z_\beta(H_k(n,p))$.

To implement this strategy, we begin by deriving upper and lower bounds
on $\Phi_{d,k}(\beta)$ via the first and the second moment method.
More precisely, in Section~\ref{Secbcrit} we are going to prove the following.

%
\begin{proposition}\label{Propbcrit}
For any $\beta$, we have
%
\begin{equation}
\label{eqPropbcrit13.1} \Phi_{d,k}(\beta)\leq\ln2+\frac{d}k\ln
\bigl(1-2^{1-k}\bigl(1-\exp(-\beta)\bigr)\bigr).
\end{equation}
Moreover, if either $d/k\leq2^{k-1}\ln2-2$ and $\beta\ge0$ or
$d/k>2^{k-1}\ln2-2$ and $\beta\le k\ln2-\ln k$, we have
%
\begin{equation}
\label{eqPropbcrit2} \Phi_{d,k}(\beta)=\ln2+\frac{d}k\ln
\bigl(1-2^{1-k}\bigl(1-\exp(-\beta)\bigr)\bigr).
\end{equation}
\end{proposition}

Since the function $\beta\in[0,\infty)\mapsto\ln2+\frac{d}k\ln
(1-2^{1-k}(1-\exp(-\beta)))$ is analytic,
it follows that the least $\beta>0$ for which~(\ref{eqPropbcrit2})
is violated marks a phase transition.
Hence, in light of (\ref{eqPropbcrit13.1}) we define
%
%
\begin{equation}
\label{eqdcrit} \beta_{\mathrm{crit}}(d,k)=\inf\biggl\{{\beta>0:
\Phi_{d,k}(\beta)<\ln2+\frac{d}k\ln\bigl(1-2^{1-k}
\bigl(1-\exp(-\beta)\bigr)\bigr)} \biggr\}.\hspace*{-30pt}
\end{equation}
We have $\beta_{\mathrm{crit}}(d,k)\in(0,\infty]$ and
Proposition~\ref{Propbcrit} readily
implies the following lower bounds on $\beta_{\mathrm{crit}}(d,k)$.

%
\begin{corollary}\label{Corbcrit}
We have $\beta_{\mathrm{crit}}(d,k)\geq k\ln2-\ln k$. If $d/k\leq
2^{k-1}\ln2-2$, then
$\beta_{\mathrm{crit}}(d,k)=\infty$.
\end{corollary}

The second main component of the proof of Theorem~\ref{Thmcond} is the
analysis of the ``cluster size'' in the planted model.
More precisely, for a hypergraph $H=(V_H,E_H)$ and a map $\sigma
:V_H\rightarrow\{{\pm1} \}$ we define
the \textit{cluster size} of $\sigma$ in $H$ as
\[
{{\mathcal C}}_\beta(H,\sigma)=\sum_{\tau\in\{{\pm1}
\}^{V_H}:\scal{\sigma
}{\tau}\geq2n/3}
\exp\bigl(-\beta E_{H}(\tau)\bigr).
\]
Thus, we sum up the contribution to the partition function of all those
maps $\tau$ whose ``overlap''
$\scal{\sigma}{\tau}=\sum_{v\in V_H}\sigma(v)\tau(v)$ with the
given $\sigma$ is big.
Concerning the cluster size in $H_k(n,p)$, there is a concentration bound
analogous to~Lemma~\ref{LemmaZAzuma}.

%
\begin{lemma}\label{LemmaMcDiarmid}
For any $\sigma: [n]\rightarrow\{{\pm1} \}$ and $\alpha
>0$, there is
$\delta=\delta(\alpha,\sigma)>0$ such that
\[
\pr\brkk{\bigl\vert\ln{{\mathcal C}}_\beta\bigl(H_k(n,p),
\sigma\bigr)-\Erw\bigl[\ln{{\mathcal C}}_\beta\bigl
(H_k(n,p),\sigma\bigr)\bigr]\bigr\vert>\alpha n}<\exp(-\delta n).
\]
\end{lemma}
\begin{pf}
This follows from McDiarmid's inequality~\cite{McDiarmid}, Theorem~3.8,
and because we have $\vert\ln{{\mathcal C}}_\beta(H,\sigma)-\ln
{{\mathcal C}}_\beta
(H',\sigma)\vert\le\beta$ for any $\sigma$ if the hypergraph $H'$~is
obtained from the hypergraph $H$ by either adding or removing a single edge.
\end{pf}

Ideally, we would like to compare the cluster size of an assignment
$\sigma$ chosen from the Boltzmann distribution on $H_k(n,p)$
with the partition function $Z_\beta(H_k(n,p))$.
Then according to the physicists' ``entropy crisis'', the condensation
phase transition should mark the point $\beta$ where
${{\mathcal C}}_\beta(H_k(n,p),\sigma)$ is of the same order of
magnitude as
$Z_\beta(H_k(n,p))$.
However, it seems difficult to calculate ${{\mathcal C}}_\beta(H_k(n,p)
,\sigma)$ directly;
the basic reason for this is that the Boltzmann distribution on a
randomly generated hypergraph is a very difficult object to approach directly.

To sidestep this difficulty, we introduce another experiment whose
outcome is much easier to study and that will emerge to be sufficient
to pin down the condensation phase transition.
This alternate experiment is the \textit{planted model}.
It is defined as follows.
Let $\bolds{\sigma}:\brk n\rightarrow\{{-1,1} \}$ be a map chosen
uniformly at random.
Moreover, given $d,k,\beta$, set
\begin{eqnarray*}
p_1 &=& \frac{\exp(-\beta)}{1-2^{1-k}(1-\exp(-\beta))}\cdot\frac
{d}{\bink{n-1}{k-1}},
\\
p_2 &=&
\frac{1}{1-2^{1-k}(1-\exp(-\beta))}\cdot\frac{d}{\bink{n-1}{k-1}}.
\end{eqnarray*}
Now, obtain a random $k$-uniform hypergraph $\G$ by inserting each
edge that is monochromatic under $\bolds{\sigma}$ with probability $p_1$
and each edge that is bichromatic under $\bolds{\sigma}$ with probability
$p_2$ independently.
In symbols, for any hypergraph $H$ with vertex set $[n]$ we have
\[
\pr\brk{\G=H\mid\bolds{\sigma}}=p_1^{E_H(\bolds{\sigma}
)}(1-p_1)^{m_1}p_2^{e(H)-E_H(\bolds{\sigma})}(1-p_2)^{m_2},
\]
where $m_1$ (resp., $m_2$) are the numbers of edges that are
monochromatic (resp., bichromatic) under $\bolds{\sigma}$ and are
\textit{not}
in $H$.

The following proposition reduces the problem of determining $\beta
_{\mathrm{crit}}(d,k)$
to that of calculating ${{\mathcal C}}_\beta(\G,\bolds{\sigma})$.
We will
prove in
Section~\ref{SecplantedCluster}.

%
\begin{proposition}\label{LemmaplantedCluster}
Assume that $ d/k = 2^{k-1}\ln2 + O_k(1)$ and $\beta_0\ge k\ln2-\ln
k$. If for all $k\ln2-\ln k\leq\beta\leq\beta_0$ we have
%
%
\begin{eqnarray}
\label{eqLemmaplantedCluster2}
\quad && \lim_{\eps\searrow0}\liminf
_{n\rightarrow\infty}\pr
\brkkk{\frac{1}n\ln{{\mathcal C}} _\beta(\G,\bolds{\sigma
})\leq
\ln2+\frac{d}k\ln\bigl(1-2^{1-k}\bigl(1-\exp(-\beta)\bigr)\bigr)-\eps}
\nonumber\\[-8pt]\\[-8pt]\nonumber
&&\qquad =1,
\end{eqnarray}
then $\beta_0\leq\beta_{\mathrm{crit}}(d,k)$.
Conversely, if
%
%
\begin{eqnarray}
\label{eqLemmaplantedCluster3}
\qquad && \lim_{\eps\searrow0}\liminf
_{n\rightarrow\infty}\pr
\brkkk{\frac{1}n\ln{{\mathcal C}}_{\beta_0}(\G,\bolds{\sigma})\geq\ln2 +
\frac{d}k\ln\bigl(1-2^{1-k}
\bigl(1-\exp(-\beta_0)\bigr)\bigr) + \eps}
\nonumber\\[-8pt]\\[-8pt]\nonumber
&&\qquad =1,
\end{eqnarray}
then $\beta_0\ge\beta_{\mathrm{crit}}(d,k)$.
\end{proposition}

Finally, in Section~\ref{SecclusterSize} we are going to estimate the
cluster size
${{\mathcal C}}_\beta(\G,\bolds{\sigma})$ to derive the following result.

%
\begin{proposition} \label{propclustersize}
Assume that $d/k =2^{k-1}\ln2 + O_k(1)$ and $\beta\geq k \ln2 - \ln k$.
Then a.a.s. the cluster size in the planted model satisfies
\[
\frac{1}{n} \ln{{\mathcal C}}_\beta(\G,\bolds{\sigma}) =
\frac{\ln
2}{2^k} - \frac
{\beta\ln2}{\exp(\beta)} + \tilde O_k
\bigl(4^{-k}\bigr).
\]
\end{proposition}

\begin{pf*}{Proof of Theorem~\ref{Thmcond}}
The result of the theorem in the case $d/k \leq2^{k-1} \ln2 -2$
follows from Corollary~\ref{Corbcrit}. Let us thus assume that $d/k =
2^{k-1} \ln2 + O_k(1)$. Because we will use Proposition~\ref
{LemmaplantedCluster}, we can also assume that $\beta\geq k \ln2-
\ln k$.
We write $c_k = d /k - 2^{k-1} \ln2+\ln2$ and $b_k = \beta- k \ln
2$. With Proposition~\ref{propclustersize}, we have a.a.s.
\begin{eqnarray*}
&&\frac{1}{n} \ln{{\mathcal C}}_\beta(\G,\bolds{\sigma}) - \biggl( \ln2 +
\frac{d}{k} \ln\bigl( 1- 2^{1-k} \bigl(1- \exp(-\beta)\bigr)\bigr)\biggr)
\\
&&\qquad = \biggl( \frac{\ln2}{2^k} - (k \ln2+b_k) \ln2 \frac{\exp(-
b_k)}{2^k}
\biggr)
\\
&&\quad\qquad{} -
\biggl(\frac{\ln2}{2^k} - \frac{c_k}{2^{k-1}} + \frac{\ln2\exp
(- b_k)}{2^k} \biggr) +
\tilde O_k\bigl(4^{-k}\bigr)
\\
&&\qquad = \frac{1}{2^{k}} \bigl[ 2 c_k - (k \ln2 + b_k +1 )
\ln2 \exp(-b_k) \bigr] + \tilde O_k\bigl(4^{-k}
\bigr)
\\
&&\qquad =  \frac{1}{2^k} \bigl[ - \Sigma_{k,d}(\beta) + \tilde
O_k\bigl(2^{-k}\bigr) \bigr].
\end{eqnarray*}
The equation $\Sigma_{k,d}(\beta)=0$ has exactly one solution $\beta
_c(d,k) \geq k \ln2- \ln k$ for $d/k >2^{k-1} \ln2 - \ln2$, and no
such solution for $d/k <2^{k-1} \ln2 - \ln2$. Moreover, $\Sigma
_{k,d}(\beta)$ is smooth for $d/k >2^{k-1} \ln2 - \ln2 + 2^{-k}$,
with derivatives of order $\Omega(k^{-4})$. Consequently, there is
$\eps_k = \tilde O_k(2^{-k})$ such that
the following is true:
\begin{longlist}[(ii)]
\item[(i)] If $d/k < 2^{k-1} \ln2- \ln2 - \eps_k$, then a.a.s.~for all
$\beta\geq k \ln2 - \ln k$,
\[
\frac{1}{n} \ln{{\mathcal C}}_\beta(\G,\bolds{\sigma}) \leq
\biggl( \ln2
+ \frac
{d}{k} \ln\bigl( 1- 2^{1-k} \bigl(1- \exp(-\beta)\bigr)
\bigr) \biggr) - {\Omega(1)}.
\]
\item[(ii)] If $d/k > 2^{k-1} \ln2- \ln2 + \eps_k $, then a.a.s.~for all
$\beta\geq k \ln2 - \ln k$:
\begin{itemize}
\item if $\beta\leq\beta_c(d,k) - \eps_k$ then
\[
\frac{1}{n} \ln{{\mathcal C}}_\beta(\G,\bolds{\sigma}) \leq
\biggl( \ln2
+ \frac
{d}{k} \ln\bigl( 1- 2^{1-k} \bigl(1- \exp(-\beta)\bigr)
\bigr) \biggr) - {\Omega(1)},
\]
\item if $\beta\geq\beta_c(d,k) + \eps_k$ then
\[
\frac{1}{n} \ln{{\mathcal C}}_\beta(\G,\bolds{\sigma}) \geq
\biggl( \ln2
+ \frac
{d}{k} \ln\bigl( 1- 2^{1-k} \bigl(1- \exp(-\beta)\bigr)
\bigr) \biggr) +{\Omega(1)}.
\]
\end{itemize}
\end{longlist}
The proof of the theorem is completed by using Proposition~\ref
{LemmaplantedCluster}.
\end{pf*}

\section{The first and the second moment}\label{Secbcrit}
\textit{Throughout this section, we assume that $0\leq d/k\leq2^{k-1}\ln
2+O_k(1)$. We let $m=\lceil dn/k\rceil$.}

In this section, we prove Proposition~\ref{Propbcrit} and also lay the
foundations for the proof of Proposition~\ref{LemmaplantedCluster}.
Recall that $H_k(n,m)$ signifies the hypergraph on $\brk n$ obtained by
choosing $m$ edges uniformly at random without replacement while for
the hypergraph $H'_k(n,m)$ we choose $m$ edges $\mathbf e_1, \ldots,
\mathbf e_m$ with replacement uniformly and independently at random,
allowing for multiple edges.

\subsection{The first moment}\label{Secfirst}
We begin with the following estimate of the first moment of $Z_\beta$
in $H'_k(n,m)$.

%
\begin{lemma}\label{LemmafirstMoment}
We have $\Erw[Z_\beta(H'_k(n,m))] = \Theta
(2^n(1-2^{1-k}(1-\exp( -\beta)))^m)$.
\end{lemma}

The proof of Lemma~\ref{LemmafirstMoment} is straightforward, but we
carry it out at leisure to introduce some notation that will be used throughout.
For a map $\sigma:\brk n\rightarrow\{{-1,1} \}$, let
\[
\Forb(\sigma)=\binkd{\bigl\vert\sigma^{-1}(-1)\bigr
\vert} {k}+\binkd{\bigl\vert\sigma^{-1}(1)\bigr\vert} {k}
\]
be the number of ``forbidden $k$-sets'' of vertices that are colored
the same under~$\sigma$.
The function $x \mapsto\bink{x}{k}+\bink{n-x}{k}$ is convex and
takes its minimal value at $x=\frac{n}2$. Therefore,
%
%
\begin{eqnarray}\label{eqForb2}
\Forb(\sigma)\ge2\binkd{n/2} {k} = 2^{1-k}N
\bigl(1+O(1/n)\bigr) =2^{1-k}N+O(N/n),
\nonumber\\[-8pt]\\[-8pt]
\eqntext{\mbox{with } N=\binkd nk.}
\end{eqnarray}

Let us call $\sigma$ \textit{balanced} if $\vert\vert\sigma
^{-1}(1)\vert-\frac{n}2 \vert\leq\sqrt n$.
Let $\Bal=\Bal_{n}$ be the set of all balanced maps $\sigma
:[n]\rightarrow
\{{\pm1} \}$. Stirling's formula yields $\vert\Bal\vert=\Omega(2^n)$.
If $\sigma\in\Bal$, then
%
%
\begin{equation}
\label{eqForb3} \Forb(\sigma)\le\binkd{n/2+\sqrt{n}} {k}+\binkd
{n/2-\sqrt{n}}
{k}=2^{1-k}N+O(N/n).
\end{equation}
For a hypergraph $H$, let
\[
\Zkb(H)=\sum_{\sigma\in\Bal}\exp\bigl(-\beta
E_{H}(\sigma)\bigr).
\]

\begin{pf*}{Proof of Lemma~\ref{LemmafirstMoment}}
By the independence of edges, we have
\begin{eqnarray*}
\Erw\brkk{\exp\bigl(-\beta E_{H'_k(n,m)}(\sigma)\bigr)}&=&\Erw\brkkkkp {\prod
_{i=1}^m\exp(-\beta\vecone_{\mathbf e_i\in\Forb(\sigma)})}
\\
&=& \prod_{i=1}^m\Erw\bigl[\exp(-\beta
\vecone_{\mathbf e_i\in\Forb
(\sigma)})\bigr]
\\
&=& \bigl({1-N^{-1}\Forb(\sigma) \bigl(1-\exp(-\beta)\bigr)}
\bigr)^m
\\
&\leq&\bigl(1-2^{1-k}\bigl(1+O(1/n)\bigr) \bigl(1-\exp(-\beta
)\bigr)
\bigr)^m.
\end{eqnarray*}
Consequently,
%
%
\begin{equation}
\label{eqLemmafirstMoment1} \Erw\bigl[Z_\beta\bigl(H'_k(n,m)
\bigr)\bigr]= O\bigl(2^n\bigl(1-2^{1-k}\bigl(1-\exp(-\beta)
\bigr)\bigr)^m\bigr).
\end{equation}
If $\sigma\in\Bal$, by (\ref{eqForb3}) we have $\Erw\brk{\exp
(-\beta E_{H'_k(n,m)}(\sigma))}=\Omega((1-2^{1-k}(1-\break \exp(-\beta)))^m)$.
Therefore,
\begin{eqnarray}
\label{eqLemmafirstMoment2}
\nonumber
\Erw\bigl[Z_\beta\bigl(H'_k(n,m)
\bigr)\bigr]&\geq&\vert\Bal\vert\cdot\Omega\bigl(\bigl(1-2^{1-k}
\bigl(1-\exp(-\beta)\bigr)\bigr)^m\bigr)
\nonumber\\[-8pt]\\[-8pt]\nonumber
&=& \Omega\bigl(2^n\bigl(1-2^{1-k}\bigl(1-\exp(-\beta)\bigr)
\bigr)^m\bigr).
\end{eqnarray}
Thus, Lemma~\ref{LemmafirstMoment} follows from~(\ref
{eqLemmafirstMoment1}) and~(\ref{eqLemmafirstMoment2}).
\end{pf*}

The following lemma relates the expectation of the partition functions
of the models $H_k(n,m)$ and $H'_k(n,m)$.

\begin{lemma}\label{Lemrelatepartition}
We have $\Erw[Z_\beta(H_k(n,m))]=\Theta(\Erw[Z_\beta(H'_k(n,m))])$.
\end{lemma}

\begin{pf}
Let $\cA$ be the event that $H'_k(n,m)$ has no multiple edges. Then,
using Fact~\ref{LemmadoubleEdge} we get
\begin{eqnarray*}
&&\Erw\bigl[Z_\beta\bigl(H'_k(n,m)\bigr)\bigr]
\ge\Erw\bigl[Z_\beta\bigl(H'_k(n,m)\bigr)\mid
\cA\bigr]\pr\brk{\cA} \ge\Erw\bigl[Z_\beta\bigl(H_k(n,m)
\bigr)\bigr]\bigl(1-o(1)\bigr),
\end{eqnarray*}
implying that
%
%
\begin{equation}
\label{eqlowerpart} \Erw\bigl[Z_\beta\bigl(H_k(n,m)\bigr)\bigr]
\le O(1)\Erw\bigl[Z_\beta\bigl(H'_k(n,m)\bigr)
\bigr].
\end{equation}
On the other hand, let $m_0=\frac{2^{1-k}\exp(-\beta
)}{1-2^{1-k}(1-\exp(-\beta))}m$ and
\[
f(x)=-x\beta-x\ln x-(1-x)\ln(1-x)+x\ln\bigl(2^{1-k}\bigr)+(1-x)\ln
\bigl(1-2^{1-k}\bigr).
\]
We observe that $f$ is strictly concave and attains its maximum at
$x=\frac{m_0}{m}$ where it is equal to $\ln(1-2^{1-k}(1-\exp(-\beta)))$.
For $\sigma\in\Bal$, we get with Stirling's formula
%
%
\begin{eqnarray}
\label{eqenergyhnm}
&&\Erw\brkk{\exp\bigl(-\beta E_{H_k(n,m)}(\sigma
)\bigr)}\nonumber
\\
\nonumber
&&\qquad =\sum_\mu\pr\brk{E_{H_k(n,m)}=\mu}
\exp(-\beta\mu)
\\
&&\qquad  \ge\sum_{\mu\in[m_0-\sqrt{m},m_0+\sqrt{m}]} \exp(-\beta\mu)
\frac{{m\choose \mu} (\Forb(\sigma))^{\mu}(N-\Forb
(\sigma))^{m-\mu}}{N^m}
\\
\nonumber
&&\qquad =\sum_{\mu\in[m_0-\sqrt{m},m_0+\sqrt{m}]} \Theta_m \biggl(
\frac{1}{\sqrt{m}} \biggr)\exp\biggl(m f \biggl(\frac
{m_0}{m} \biggr)
\biggr)\Theta(1)
\\
\nonumber
&&\qquad =\Theta \bigl(1-2^{1-k}\bigl(1-\exp(-\beta)
\bigr)^m\bigr).
\end{eqnarray}
Therefore,
%
%
\begin{eqnarray}
\label{equpperpart} \Erw\bigl[Z_\beta\bigl(H_k(n,m)\bigr)\bigr]
&\ge&\vert\Bal\vert\cdot\Erw\brkk{\exp\bigl(-\beta
E_{H_k(n,m)}(\sigma)
\bigr)}
\nonumber\\[-8pt]\\[-8pt]\nonumber
& =& \Omega\bigl(2^n \bigl(1-2^{1-k}\bigl(1-\exp(-
\beta)\bigr)^m\bigr)\bigr).
\end{eqnarray}
Combining (\ref{eqlowerpart}), Lemma~\ref{LemmafirstMoment} and
(\ref{equpperpart}) proves the assertion.
\end{pf}

As a further consequence of Lemma~\ref{LemmafirstMoment}, we obtain
the following.

%
\begin{corollary}\label{CorfirstMoment}
1. We have $\Phi_{d,k}(\beta)\leq
\ln2+\frac{d}k\ln(1-2^{1-k}(1-\exp(-\beta)))$ for all $d,\beta$.
2. Assume that $d,\beta$ are such that
\[
\limsup_{n\rightarrow\infty}\frac{1}n\Erw\bigl[\ln\Zkc
\bigl(H'_k(n,m)\bigr)\bigr]<\ln2+\frac{d}k\ln
\bigl(1-2^{1-k}\bigl(1-\exp(-\beta)\bigr)\bigr).
\]
Then
$\Phi_{d,k}(\beta)<\ln2+\frac{d}k\ln(1-2^{1-k}(1-\exp(-\beta)))$.
\end{corollary}

\begin{pf}
Let ${{\mathcal E}}$ be the event that $\vert e(H_k(n,p))-m\vert\leq
\sqrt n\ln n$.
Then we can couple the random hypergraphs $H_k(n,m)$ and $H_k(n,p)$ given
${{\mathcal E}}$ as follows.
\begin{enumerate}
\item Choose a random hypergraph $H_0=H_k(n,m)$.
\item Let $e=\Bin(\bink nk,p )$ be a binomial random
variable given that $\vert e-m\vert\leq\sqrt n\ln n$.
\item Obtain a random hypergraph $H_1$ from $H_0$ as follows:
\begin{itemize}
\item If $e\geq m$, choose a set of $e-m$ random edges from all edges
not present in $H_0$ and add them to $H_0$.
\item If $e<m$, remove $m-e$ randomly chosen edges from $H_0$.
\end{itemize}
\end{enumerate}
The\vspace*{1pt} outcome $H_1$ has the same distribution as $H_k(n,p)$ given
${{\mathcal E}}$,
and $H_0,H_1$ differ in at most $\sqrt n\ln n$ edges.
Therefore, noting that $\frac{1}n\vert\ln Z_\beta\vert\leq\frac
{d}k\beta+\ln
2$ with certainty, we obtain with Fact~\ref{eqLip}:
%
%
\begin{eqnarray}\label{eqJensen1}
\frac{1}n\Erw\ln Z_\beta\bigl(H_k(n,p)\bigr)&\leq&
\frac{1}n\Erw\bigl[\ln Z_\beta(H_1)\bigr]+ \biggl(
\frac{d}k\beta+\ln2 \biggr)\pr\brk{\neg{{\mathcal E}}}
\nonumber
\\
&\leq&\frac{1}n\Erw\bigl[\ln Z_\beta(H_0)\bigr]+
\frac{\beta\ln n}{\sqrt
n}+ \biggl(\frac{d}k\beta+\ln2 \biggr)\pr\brk{\neg{{
\mathcal E}}}
\\
&=&\frac{1}n\Erw\bigl[\ln Z_\beta\bigl(H_k(n,m)
\bigr)\bigr]+ \biggl(\frac{d}k\beta+\ln2 \biggr)\pr\brk{\neg
{{\mathcal
E}}}+o(1).\nonumber
\end{eqnarray}
Since $e(H_k(n,p))$ is a binomial random variable with mean $m+O(1)$,
Lemma~\ref{lemchernoffsum} implies that $\pr\brk{\neg
{{\mathcal E}}}=o(1)$.
Thus, by~(\ref{eqJensen1}) and Jensen's inequality,
\begin{eqnarray*}
\frac{1}n\Erw\ln Z_\beta\bigl(H_k(n,p)\bigr) &\le&
\frac{1}n\Erw\bigl[\ln Z_\beta\bigl(H_k(n,m) \bigr)
\bigr]+o(1)
\\
&\le&\frac{1}n\ln\Erw\bigl[Z_\beta\bigl(H_k(n,m)
\bigr)\bigr]+o(1).
\end{eqnarray*}
The first assertion follows by Lemmas~\ref{LemmafirstMoment} and~\ref{Lemrelatepartition} and taking $n\rightarrow\infty$. Also
the second
assertion readily follows.
\end{pf}

We conclude this section by observing that the contribution to $Z_\beta
$ of certain ``exotic'' $\sigma$ is negligible.
We begin with $\sigma$ that are very imbalanced.

%
\begin{lemma}\label{Lemimbalanced}
For any $\eps>0$ there is $\delta>0$ such that the following is true.
Let $\bar B_\eps$ be the set of all $\sigma:\brk n\rightarrow\{
{\pm1} \}$
such that $\vert\vert\sigma^{-1}(1)\vert-\frac{n}2\vert>\eps n$.
Moreover, let
\[
Z_{\beta,\bar B_\eps}(H)=\sum_{\sigma\in\bar B_\eps}\exp\bigl
(-\beta
E_H(\sigma)\bigr).
\]
Then $\Erw[Z_{\beta,\bar B_\eps}(H_k(n,m))]\leq\exp(-\delta n)\Erw
[Z_{\beta}(H_k(n,m))]$.
\end{lemma}

\begin{pf}
Stirling's formula implies that for any $\eps>0$ there is $\delta>0$
such that $\frac{1}n\ln\vert\bar B_\eps\vert<\ln2-\delta$.
Hence, (\ref{eqForb2}) implies together with the independence of the
edges that
\begin{eqnarray*}
\Erw\bigl[Z_{\beta,\bar B_\eps}\bigl(H'_k(n,m)\bigr)
\bigr]&=&\sum_{\sigma\in\bar
B_\eps
}\Erw\bigl[\exp\bigl(-\beta
E_{H'_k(n,m)}(\sigma)\bigr)\bigr]
\\
&\leq&\vert\bar B_\eps\vert\bigl(1-2^{1-k}\bigl(1-
\exp(-\beta)\bigr)\bigr)^m
\\
&\leq&\exp(-\delta n)2^n\bigl(1-2^{1-k}\bigl(1-\exp(-\beta)
\bigr)\bigr)^m.
\end{eqnarray*}
The assertion follows from the remark that [as in equation (\ref
{eqlowerpart})]
\[
\Erw\bigl[Z_{\beta,\bar B_\eps}\bigl(H_k(n,m)\bigr)\bigr] =
O\bigl(\Erw\bigl[Z_{\beta,\bar B_\eps}\bigl(H'_k(n,m)\bigr)\bigr]\bigr),
\]
and from Lemma~\ref
{Lemrelatepartition}.
\end{pf}

%
\begin{lemma}\label{LemavgEngergy}
For any $\eps>0$, there is $\delta>0$ such that the following is true.
Let $m_0=\frac{2^{1-k}\exp(-\beta)}{1-2^{1-k}(1-\exp(-\beta))}m$ and
\begin{eqnarray*}
Z_{\beta,\eps}(H)&=&\sum_{\sigma:\brk n\rightarrow\{{\pm
1} \}}\exp\bigl(-\beta
E_H(\sigma)\bigr)\cdot\vecone_{\vert E_H(\sigma)-m_0\vert>\eps m}.
\end{eqnarray*}
Then $\Erw[Z_{\beta,\eps}(H_k(n,m))]\leq\exp(-\delta n)\Erw
[Z_{\beta
}(H_k(n,m))]$.
\end{lemma}

\begin{pf}
Let $M_0= \{{\mu\in\brk m:\vert\mu-m_0\vert>\eps m} \}$.
Moreover, for
$\alpha>0$ let $B_{\alpha}$ be the set of all
$\sigma:\brk n\rightarrow\{{\pm1} \}$ such that $\vert
\vert\sigma
^{-1}(1)\vert-\frac{n}2 \vert<\alpha n$.
Then by Lemma~\ref{Lemimbalanced} there exists $\delta>0$ such that
%
%
\begin{eqnarray}
\label{eqavgEng1} \Erw\bigl[Z_{\beta,\eps}\bigl(H_k(n,m)\bigr
)\bigr]&
\leq&\exp(-\delta n)\Erw\bigl[Z_{\beta
}\bigl(H_k(n,m)\bigr)
\bigr]
\nonumber\\[-8pt]\\[-8pt]\nonumber
&&{} +\sum_{\mu\in M_0}\sum
_{\sigma\in B_\alpha}\exp(-\beta\mu)\pr\brkk{E_{H_k(n,m)}(\sigma
)=\mu}.
\end{eqnarray}
As in the proof of Lemma~\ref{Lemrelatepartition}, we define
$f(x)=-x\beta-x\ln x-(1-x)\ln(1-x)+x\ln(2^{1-k})+(1-x)\ln
(1-2^{1-k})$ and find that for any $\gamma>0$ we can choose $\alpha
>0$ small enough so that
\[
\frac{1}m\ln\bigl(\exp(-\beta\mu)\pr\brkk{E_{H_k(n,m)}(\sigma)=
\mu} \bigr)\leq\gamma+f \biggl(\frac{\mu}{m} \biggr)\qquad
\mbox{for all }
\sigma\in B_\alpha.
\]
Because $f$ is strictly concave and attains its maximum at $x=\frac
{m_0}{m}$, there is $\delta'>0$ such that
%
%
\begin{equation}
\label{eqavgEng2} \sum_{\mu\in M_0}\sum
_{\sigma\in B_\alpha}\exp(-\beta\mu)\pr\bigl[E_{H_k(n,m)}(\sigma
)=\mu
\bigr]\leq\exp\bigl(-\delta' n\bigr)\Erw\bigl[Z_{\beta}
\bigl(H_k(n,m)\bigr)\bigr].\hspace*{-25pt}
\end{equation}
Finally, the assertion follows from~(\ref{eqavgEng1}) and~(\ref{eqavgEng2}).
\end{pf}

\subsection{The second moment}\label{Secsecond}

In Section~\ref{Secfirst}, we derived an upper bound on $\Phi
_{d,k}(\beta)$ by calculating the expectation of $Z_{\beta}(H'_k(n,m))$
(cf. Corollary~\ref{CorfirstMoment}).
Here, we obtain for certain values of $\beta$ and $d$ a matching lower
bound by estimating the second moment $\Erw[\Zkb(H'_k(n,m))^2]$.
To this end, we define for $\alpha\in[-1,1]$,
%
%
\begin{equation}
\label{eqZbetaalpha} Z_{\beta}(\alpha)=\sum_{\sigma,\tau\in\Bal
:\scal{\sigma}{\tau
}=\alpha n}
\exp\bigl(-\beta\bigl(E_{H'_k(n,m)}(\sigma)+E_{H'_k(n,m)}(\tau
)\bigr)
\bigr).
\end{equation}
Thus, in~(\ref{eqZbetaalpha}) we sum over balanced pairs $\sigma,\tau
:\brk n\rightarrow\{{\pm1} \}$
that agree on precisely $n((1+\alpha)/2)$ vertices. Hence, we can
express the second moment as
\begin{eqnarray*}
\Erw\brkk{\Zkb\bigl(H'_k(n,m)\bigr)^2}&=&
\sum_{\sigma,\tau\in\Bal}\Erw\bigl[\exp\bigl(-\beta
\bigl(E_{H'_k(n,m)}(\sigma)+E_{H'_k(n,m)}(\tau)\bigr)\bigr)\bigr]
\\
&=&\sum_{\nu=0}^n\Erw\bigl[Z_\beta(2
\nu/n-1)\bigr].
\end{eqnarray*}
Consequently, we need to bound $Z_{\beta}(\alpha)$ for $-1\leq\alpha
\leq1$.
Recall that $\cH(z)=-z\ln z-(1-z)\ln(1-z)$.

%
\begin{lemma}\label{LemmaLambda}
For $\alpha\in[-1,1]$, we have
\begin{eqnarray*}
\frac{1}n{\ln\Erw\bigl[Z_{\beta}(\alpha)\bigr]}&=&\ln2+
\Lambda_\beta(\alpha)-\frac{\ln n}{2n}+O(1/n),
\end{eqnarray*}
where
\begin{eqnarray*}
\Lambda_{\beta}(\alpha)& =& \cH\biggl({\frac{1+\alpha}{2}} \biggr)
+ \frac{d}k\ln\biggl[1-2^{1-k}\bigl(1-\exp(-\beta)\bigr)
\\
&&{}\times \brkkk{2-
\bigl(1-\exp(-\beta)\bigr)\frac{(1+\alpha)^k+(1-\alpha
)^k}{2^k}} \biggr].
\end{eqnarray*}
\end{lemma}

\begin{pf}
Let $e$ be a randomly chosen edge.
Let $\sigma,\tau:\brk n\rightarrow\{{\pm1} \}$ be two
balanced maps with
overlap $\scal\sigma\tau=\alpha n$.
Let us write $\sigma\vDash e$ if $e\notin\Forb(\sigma)$ (i.e., $e$
is bichromatic under $\sigma$).
By inclusion--exclusion,
\begin{eqnarray*}
\pr\brkk{\sigma\vDash e}, \pr\brkk{\tau\vDash e}&=&1-2^{1-k} +O(1/n),
\\
\pr\brkk{\sigma, \tau\vDash e}&=&1-2^{2-k}+2^{1-2k} \bigl((1+
\alpha)^k+ (1-\alpha)^k \bigr)+O(1/n).
\end{eqnarray*}
Hence, by the independence of edges,
%
%
\begin{eqnarray}\label{eqLemmaLambda1}
\Erw\brkk{Z_{\beta}(\alpha)}&=&\sum_{\sigma,\tau:\scal\sigma\tau
=\alpha n}\Erw
\prod_{i=1}^m\exp\brkk{-\beta(
\vecone_{\sigma
\nvDash e_i}+\vecone_{\tau\nvDash e_i})}
\nonumber
\\
&=& \sum_{\sigma,\tau:\scal\sigma\tau=\alpha n} \bigl({\Erw\brkk
{\exp\brkk{-\beta(
\vecone_{\sigma\nvDash e_1}+\vecone_{\tau
\nvDash e_1})}}} \bigr)^m
\nonumber
\\
&=&2^n\binkd{n} {(1+\alpha)n/2}\bigl(\pr\brkk{\sigma,\tau\vDash
e_1}\nonumber
\\
&&{}+\exp(-\beta) \bigl(\pr\brkk{\sigma\vDash e_1,\tau
\nvDash e_1}
+\pr\brkk{\sigma\nvDash e_1,\tau\vDash e_1}\bigr)
\\
&&{} +\exp(-2
\beta)\cdot\pr\brkk{\sigma,\tau\nvDash e_1}\bigr)^m\nonumber
\\
&=&2^n\binkd{n} {(1+\alpha)n/2}\bigl(1+O(1/n)\bigr)
\bigl[1-2^{2-k}\bigl(1-\exp(-\beta)\bigr)\nonumber
\\
&&{}+2^{1-2k}
\bigl(1-\exp(-\beta)\bigr)^2\bigl((1+\alpha)^k+(1-\alpha
)^k\bigr) \bigr]^m.\nonumber
\end{eqnarray}
Furthermore, by Stirling's formula,
%
%
\begin{equation}
\label{eqLemmaLambda2} \binkd{n} {(1+\alpha)n/2}=O\bigl
(n^{-1/2}\bigr)\exp
\biggl(n\cH\biggl(\frac
{1+\alpha}{2} \biggr) \biggr).
\end{equation}
The assertion follows by combining (\ref{eqLemmaLambda1}) and~(\ref
{eqLemmaLambda2}).
\end{pf}

Hence, we need to study the function $\Lambda_\beta$.
Since $\Lambda_\beta(\alpha)=\Lambda_\beta(-\alpha)$, $\alpha=0$
is a stationary point.
Moreover, with
\[
s = s(\alpha,\beta) =
1-2^{1-k}\bigl(1- \exp(-\beta)\bigr) \biggl[ 2 - \bigl(1-\exp
(-\beta)
\bigr) \frac{(1+\alpha)^k + (1-\alpha)^k}{2^k} \biggr]
\]
the first two derivatives of $\Lambda_\beta$ work out to be
%
%
\begin{eqnarray}
\Lambda_\beta'(\alpha)&=&\frac{\ln(1-\alpha)-\ln(1+\alpha
)}{2}
\nonumber\\[-8pt]\label{eqLambda'} \\[-8pt]
&&{} +\frac{2d}{4^k s} {\bigl(\exp(-\beta)-1\bigr)^2\bigl((1+
\alpha)^{k-1}-(1-\alpha)^{k-1}\bigr)},\nonumber
\\
\Lambda_\beta''(\alpha)&=&
\frac{1}{\alpha^2-1}+\frac{2d(k-1)(\exp
(-\beta)-1)^2}{4^ks} \bigl({(1+\alpha)^{k-2}+(1-
\alpha)^{k-2}} \bigr)\hspace*{-30pt}
\nonumber\\[-8pt]\label{eqLambda''} \\[-8pt]\nonumber
&&{}-\frac{dk(1-\exp(-\beta))^4}{2^{4k-2}s^2}{ \bigl[(1+\alpha)^{k-1}-(1-
\alpha)^{k-1} \bigr]^{2} }.
\end{eqnarray}
In particular,
%
%
\begin{equation}
\label{eqlocalMax} \Lambda_\beta''(0)=-1+\tilde
O_k\bigl(2^{-k}\bigr)<0.
\end{equation}
Hence, there is a local maximum at $\alpha=0$. As a consequence, we
have
\[
\Erw\bigl[\Zkc\bigl(H'_k(n,m)\bigr)^2\bigr]= O\bigl(\Erw\bigl[\Zkc\bigl(H'_k(n,m)\bigr)\bigr]^2\bigr),
\]
if
$\Lambda
_\beta$ has a strict \textit{global} maximum at $\alpha=0$.
More generally, we have
the following.

\begin{lemma}\label{LemmaLaplace}
Assume that $\beta\ge0$ and $J\subset\brk{-1,1}$ is a compact set
such that
$\Lambda_\beta(\alpha)<\Lambda_\beta(0)$ for all $\alpha\in
J\setminus\{{0} \}$.
Then
\[
\sum_{\nu=0}^n\Erw\brkk{\Zkc(2\nu/n-1)}
\vecone_{2\nu/n-1\in J}= O\bigl(\Erw\bigl[\Zkc\bigl(H'_k(n,m)
\bigr)\bigr]^2\bigr).
\]
\end{lemma}

\begin{pf}
We start by observing that
$\frac{\ln2+\Lambda_\beta(0)}2=\ln2+\frac{d}k\ln
(1-2^{1-k}(1-\exp
(-\beta)))$.
Hence, Lemma~\ref{LemmafirstMoment} yields
%
%
\begin{equation}
\label{eqLemmaLaplace0} \exp\brkk{n\bigl(\ln2+\Lambda_\beta
(0)\bigr)}=O\bigl(\Erw
\bigl[\Zkc\bigl(H'_k(n,m)\bigr)\bigr]^2
\bigr).
\end{equation}
Now, by~(\ref{eqlocalMax}),
there exist $\eta,c>0$ such that
$\Lambda_\beta(\alpha)\leq\Lambda_\beta(0)-c\alpha^2\mbox{ for
all }\alpha\in J_0=J\cap(-\eta,\eta)$.
Hence, by Lemma~\ref{LemmaLambda} and~(\ref{eqLemmaLaplace0})
\begin{eqnarray}
\label{eqLemmaLaplace2} && \sum_{\nu=0}^n\Erw\brkk{
\Zkc(2\nu/n-1)}\vecone_{2\nu/n-1\in
J_0}\nonumber
\\
&&\qquad = O\bigl(n^{-1/2}2^n\bigr)\sum_{\nu=0}^n
\exp\bigl(n\Lambda_\beta(2\nu/n-1)\bigr)\vecone_{2\nu/n-1\in J_0}
\nonumber\\[-8pt]\\[-8pt]\nonumber
&&\qquad = O\bigl(2^n\exp\bigl(n\Lambda_\beta(0)\bigr)\bigr)\sum
_{\nu:\vert2\nu/n-1\vert<\eta}\frac
{\exp(-nc(2\nu/n-1)^2)}{\sqrt n}
\nonumber
\\
&&\qquad =O\bigl(2^n\exp\bigl(n\Lambda_\beta(0)\bigr)\bigr)= O
\bigl(\Erw\bigl[\Zkc\bigl(H'_k(n,m)\bigr)
\bigr]^2\bigr).
\nonumber
\end{eqnarray}
Further, let $J_1=J\setminus(-\eta,\eta)$.
Then $J_1$ is compact.
Hence, there exists $\delta>0$ such that $\Lambda_\beta(\alpha
)<\Lambda_\beta(0)-\delta$ for all $\alpha\in J_1$.
Therefore, Lemma~\ref{LemmaLambda} and~(\ref{eqLemmaLaplace0}) yield
%
%
\begin{eqnarray}
\label{eqLemmaLaplace3} \sum_{\nu=0}^n\Erw\brkk{
\Zkc(2\nu/n-1)}\vecone_{2\nu/n-1\in J_1} &=&O\bigl(n2^n\bigr)\sup
_{\alpha\in J_1}\exp\bigl(n\Lambda_\beta(\alpha)\bigr)\nonumber
\\
&=&O\bigl(n2^n\bigr)\exp\bigl(n\bigl(\Lambda_\beta(0)-
\delta\bigr)\bigr)
\\
&=&O\bigl(\Erw\bigl[\Zkc\bigl(H'_k(n,m)\bigr)
\bigr]^2\bigr).
\nonumber
\end{eqnarray}
Finally, the assertion follows from~(\ref{eqLemmaLaplace2}) and~(\ref
{eqLemmaLaplace3}).
\end{pf}
Now we prove that for the set $J$ from Lemma~\ref{LemmaLaplace} we
have at least $[-1+2^{-3k/4},1-2^{-3k/4}] \subset J$ for all $\beta\ge0$.

\begin{lemma}\label{Lemmasmm1}
For $d/k=2^{k-1}\ln2+O_k(1)$ and $\beta\ge0$ we have $\Lambda_\beta
(\alpha)<\Lambda_\beta(0)$ for all $\alpha\ne0$ with $\vert
\alpha\vert\le
1-2^{-3k/4}$.
\end{lemma}

\begin{pf}
We know that there is a local maximum at $\alpha=0$.
Moreover, we read off of~(\ref{eqLambda''}) that $\Lambda_\beta
''(\alpha)<0$ if $\vert\alpha\vert<1-6\ln k/k$, and thus
\[
\Lambda_\beta(0)>\Lambda_\beta(\alpha)\qquad\mbox{for all }
\alpha\in\bigl(-(1-6\ln k/k),1-6\ln k/k\bigr).
\]
Further, we obtain from~(\ref{eqLambda'}) for $\vert\alpha\vert\ge
1-6\ln k/k$
\begin{eqnarray*}
\Lambda_\beta'(\alpha)&\leq&\frac{\ln(1-\alpha)}2+
\frac
{2d(1-\exp(-\beta))^2(1+\alpha)^{k-1}}{4^k(1+O_k(2^{-k}))}
\\
&\leq&\frac{\ln(1-\alpha)}2+\frac{d(1-\exp(-\beta))^2\exp
((1+\alpha)(k-1)/2)}{2^k(1+O_k(2^{-k}))}.
\end{eqnarray*}
Hence, for $k$ large enough $\Lambda_\beta'(\alpha)<0$ if $\vert
\alpha
\vert<1-2.01\ln k/k$ and a similar estimate yields
%
%
\begin{equation}
\label{eqder} \Lambda_\beta'(\alpha)>0\qquad\mbox{if }
\vert\alpha\vert>1-1.99\ln k/k.
\end{equation}

Thus, to proceed we need to evaluate $\Lambda_\beta$ at $\vert
\alpha
\vert=1-\gamma\ln k/k$ for $\gamma\in[1.99,2.01]$ and at $\vert
\alpha\vert=1-2^{-3k/4}$.
We find
\[
\Lambda_\beta(\alpha)=-\ln2+o_k(1)
\]
for $\vert\alpha\vert=1-\gamma\ln k/k$ with $\gamma\in
[1.99,2.01]$ and
$\Lambda_\beta(\alpha)=-\ln2+o_k(1)$
for $\vert\alpha\vert=1-2^{-3k/4}$ proving the assertion.
\end{pf}

%
\begin{lemma}\label{LemmaFeli}
The function $\beta\mapsto\Lambda_\beta(\alpha)-\Lambda_\beta
(0)$ is nondecreasing for $\alpha\ne0$. In particular, if $d>0$ and
$\beta_0\ge0$ are such that
$\Lambda_{\beta_0}(\alpha)<\Lambda_{\beta_0}(0)$ for all $\alpha
\neq0$, then $\Lambda_\beta(\alpha)<\Lambda_\beta(0)$ for all
$\alpha\neq0,0\le\beta<\beta_0$.
\end{lemma}

\begin{pf}
The derivative of $\Lambda_\beta$ with respect to $\beta$ works out
to be
\begin{eqnarray*}
\hspace*{-5pt}&& \frac{\partial\Lambda_\beta}{\partial\beta}
\\
\hspace*{-5pt}&&\hspace*{4pt}\quad = \frac{d}k
\cdot \frac{2^{2-2k}((1 + \alpha)^k + (1 - \alpha)^k)\exp(-\beta)(1 - \exp
(-\beta))-2^{2-k}\exp(-\beta)}{
1 - 2^{2-k}(1 - \exp
(-\beta))+2^{1-2k}(1 - \exp(-\beta))^2((1 + \alpha)^k + (1 - \alpha)^k)}.
\end{eqnarray*}
Substituting $z=(1+\alpha)^k+(1-\alpha)^k$ and $b=1-\exp(-\beta)$
in the above, we obtain
\[
g(z)=\frac{d}k\cdot\frac{2^{2-2k}b(1-b)z-2^{2-k}(1-b)}{
1-2^{2-k}b+2^{1-2k}b^2z}.
\]
Because a function $z \mapsto\frac{az-b}{cz+d}$ with $a,b,c,d \ge0$
is nondecreasing, this completes the proof.
\end{pf}
With these instruments in hand we identify regimes of $d$ and $\beta$
where $\Lambda_\beta(\alpha)$ takes its global maximum at $\alpha=0$.

%
\begin{lemma}\label{Lemmasmm2}
Assume that $d/k=2^{k-1}\ln2+O_k(1)$ and $\beta\le k\ln2- \ln k$.
Then $\Lambda_\beta(0)>\Lambda_\beta(\alpha)$ for all $\alpha\in
[-1,1]\setminus\{{0} \}$.
\end{lemma}

\begin{pf}
For $\vert\alpha\vert\le1-2^{-3k/4}$ this is the statement of
Lemma~\ref
{Lemmasmm1}. We write $\alpha=1-\delta$ with $\delta\in[0,2^{-3k/4}]$.
Let
\[
f_\beta(\delta) = \bigl(1 - \exp(-\beta) \bigr) \biggl[ 2- \bigl(1-
\exp(-\beta) \bigr) \frac{(2-\delta)^k+\delta^k}{2^k} \biggr]
\in[0,2].
\]
For $\beta= k \ln2 - \ln k$, we have the expansion
\begin{eqnarray*}
f_\beta(\delta) &=& \biggl(1 - \frac{k}{2^k} \biggr) \biggl[ 2-
\biggl(1- \frac{k}{2^k} \biggr) \biggl( 1-k \frac{\delta}{2} +
\tilde
O_k\bigl(4^{-k}\bigr) \biggr) \biggr]
\\
&=& 1 + k \frac{\delta}{2} + \tilde O_k\bigl(4^{-k}
\bigr).
\end{eqnarray*}

Therefore,
\begin{eqnarray*}
\Lambda_\beta(\alpha)& =& - \frac{\delta}{2} \ln\biggl(
\frac
{\delta}{2} \biggr) - \biggl(1-\frac{\delta}{2} \biggr) \ln
\biggl(1-
\frac{\delta}{2} \biggr)
\\
&&{} + \frac{d}k \ln\brkkk{1- 2^{1-k} \brkkk{1 + k
\frac{\delta}{2} + \tilde O_k\bigl(4^{-k}\bigr)}}
\\
&=& - \ln2 - \frac{\delta}{2} \ln\delta+ \frac{ \delta}{2} - (k-1)
\frac{\delta}{2}\ln2 + O_k\bigl(2^{-k}\bigr).
\end{eqnarray*}

The function $\delta\mapsto- \frac{\delta}{2} \ln\delta+ \frac{
\delta}{2} - (k-1) \frac{\delta}{2}\ln2$ is easily studied: it
takes its maximum at $\delta_0 = 2^{1-k}$ for which it is equal to
$2^{-k}$. Hence, for $\alpha= 1-\delta$ with $\delta\in[0,2^{-3k/4}]$,
\[
\Lambda_\beta(\alpha) \leq- \ln2 + O_k
\bigl(2^{-k}\bigr). %
\]
By symmetry, this also holds for $\alpha= -1+\delta$ with $\delta\in
[0,2^{-3k/4}]$.
By comparison,
\begin{eqnarray*}
\Lambda_\beta(0)&=&\ln2+ \bigl( 2^{k-1} \ln2 +
O_k(1) \bigr) \ln\biggl({1-2^{2-k}+\frac{4k}{4^k}+
O_k\bigl(4^{-k}\bigr)} \biggr)
\\
&=&-\ln2+ 2^{1-k} k\ln2 +O_k\bigl(2^{-k}\bigr).
\end{eqnarray*}
Therefore, $\Lambda_\beta(0)>\Lambda_\beta(\alpha)$ for all
$\alpha\neq0$ if $\beta= k\ln2-\ln k$. Using Lemma~\ref
{LemmaFeli}, we can expand the result to all $\beta\leq k\ln2-\ln k$.
\end{pf}

%
\begin{lemma}\label{LemmazeroTermpSmm}
Assume that $d/k \le2^{k-1}\ln2-2$ and $\beta\ge0$.
Then $\Lambda_\beta(0)>\Lambda_\beta(\alpha)$ for all $\alpha\in
[-1,1]\setminus\{{0} \}$.
\end{lemma}
\begin{pf}
Let $r_k = O_k(1)$ such that $d/k = 2^{k-1}\ln2 + r_k$. Define the
function $\Lambda_\infty:[-1,1]\rightarrow\RR$ as
\[
\alpha\mapsto\cH\biggl(\frac{1+\alpha}{2} \biggr)+\frac{d}k\ln
\bigl({1-2^{2-k}+2^{1-2k}\bigl((1+\alpha)^k+(1-
\alpha)^k\bigr)} \bigr).
\]
Analogously to the proof of Lemma~\ref{Lemmasmm2}, we get
$ \Lambda_\infty(\alpha) \leq- \ln2 - ( \ln2 + 2 r_k-1) 2^{-k} +
\tilde O_k(4^{-k})$ for all $\alpha$
and
$\Lambda_\infty(0)=-\ln2- 2( \ln2 + 2 r_k) 2^{-k} + \tilde O_k(4^{-k})$,
which implies that for $r_k \le-2$ we have
$\Lambda_\infty(\alpha)<\Lambda_\infty(0)$ for all $\alpha\in
[-1,1]\setminus\{{0} \}$.
Because the continuous functions $\Lambda_\beta$ converge uniformly
to $\Lambda_\infty$ as $\beta\rightarrow\infty$,
we conclude that there is $\beta_0\ge0$ such that for all $\beta
>\beta_0$,
%
%
\begin{equation}
\label{eqFeli} \Lambda_\beta(\alpha)<\Lambda_\beta(0)\qquad
\mbox{for all }\alpha\in[-1,1]\setminus\{{0} \}.
\end{equation}
Hence, Lemma~\ref{LemmaFeli} implies that~(\ref{eqFeli}) holds for
all $\beta\ge0$, as desired.
\end{pf}

\begin{pf*}{Proof of Proposition~\ref{Propbcrit}}
The first assertion follows directly from Corollary~\ref{CorfirstMoment}.
Moreover, if $d,\beta$ are such that for some $n$-independent number
$C>0$ we have
%
%
\begin{equation}
\label{eqPropbcrit1} \Erw\bigl[\Zkc\bigl(H'_k(n,m)
\bigr)^2\bigr]\leq C\cdot\Erw\bigl[\Zkc\bigl(H'_k(n,m)
\bigr)\bigr]^2,
\end{equation}
then the Paley--Zygmund inequality implies that
%
%
\begin{eqnarray}\label{eqPZ}
\pr\brkk{\Zkc\bigl(H'_k(n,m)\bigr)\geq
\Erw\bigl[\Zkc\bigl(H'_k(n,m)\bigr)\bigr]/2}&\geq&
\frac
{\Erw
[\Zkc(H'_k(n,m))]^2}{4\Erw[\Zkc(H'_k(n,m))^2]}
\nonumber\\[-8pt]\\[-8pt]\nonumber
&\geq&\frac{1}{4C} >0.
\end{eqnarray}
Let $\cA$ be the event that $H'_k(n,m)$ has no multiple edges. Since
$\cA
$ occurs a.a.s. by Fact~\ref{LemmadoubleEdge}, (\ref{eqPZ}) implies that
%
%
\begin{equation}
\label{eqPropbcrit3} \pr\brkk{\Zkc\bigl(H'_k(n,m)\bigr)\geq
\Erw\bigl[\Zkc\bigl(H'_k(n,m)\bigr)\bigr]/2\mid\cA}\geq
\frac{1-o(1)}{4C}.
\end{equation}
Further, since the number $e(H_k(n,p))$ of edges in $H_k(n,p)$ has a binomial
distribution with mean $m+O(1)$, Stirling's formula implies that\break
$\pr\brk{e(H_k(n,p))= m}\geq\Omega(n^{-1/2})$.
Because given $e(H_k(n,p))=m$, $H_k(n,p)$ is identically distributed as
$H'_k(n,m)$ given $\cA$, (\ref{eqPropbcrit3}) implies that
%
%
\begin{equation}
\label{eqPropbcrit4} \pr\brkk{\Zkc\bigl(H_k(n,p)\bigr)\geq\Erw
\bigl[\Zkc
\bigl(H'_k(n,m)\bigr)\bigr]/2}\geq\Omega
\bigl(n^{-1/2}\bigr).
\end{equation}
The concentration bound from Lemma~\ref{LemmaZAzuma} and (\ref
{eqPropbcrit4}) yields $\ln\Erw[\Zkc(H'_k(n,m))]-\Erw
[\ln
\Zkc(H_k(n,p))]-\ln2=o(n)$.
Hence, if (\ref{eqPropbcrit1}) is true, then
%
%
\begin{equation}
\label{eqPropbcrit5} \frac{1}n\Erw\brkk{\ln\Zkc\bigl
(H_k(n,p)\bigr)}
\geq\frac{1}n\ln\Erw\bigl[\Zkc\bigl(H'_k(n,m)
\bigr)\bigr]-o(1).
\end{equation}

Finally, Lemma~\ref{LemmaLaplace} and Lemma~\ref{LemmazeroTermpSmm}
imply that~(\ref{eqPropbcrit1}) holds
for all $\beta\ge0$ and $d/k \le2^{k-1}\ln2-2$.
Moreover, by Lemma~\ref{LemmaLaplace} and Lemma~\ref{Lemmasmm2} the
bound (\ref{eqPropbcrit1}) is true
if $d/k= 2^{k-1}\ln2+O_k(1)$ and $\beta\leq k\ln2-\ln k$.
Thus, the assertion follows from (\ref{eqPropbcrit5}).
\end{pf*}

\section{The planted model}\label{SecplantedCluster}

\textit{The aim of this section is to prove Proposition}~\ref
{LemmaplantedCluster}. \textit{Throughout the section}, \textit{we let} $m=\lceil
dn/k\rceil$. \textit{For} $\eps>0$, \textit{we let} $B_\eps$ \textit{be the set of all} $\sigma
:\brk n\rightarrow\{{\pm1} \}$ \textit{such that} $\vert\vert\sigma
^{-1}(1)\vert-\frac{n}2
\vert<\eps n$. \textit{Further}, \textit{we let} $\bolds{\sigma}:\brk n\rightarrow\{
{\pm1} \}$ \textit{be a
map chosen uniformly at random and $\G$ be the random hypergraph
obtained by inserting each edge that is monochromatic under} $\bolds
{\sigma}$
\textit{with probability} $p_1$ \textit{and each edge that is bichromatic with
probability} $p_2$.

\subsection{Quiet planting}

We begin with the second part of Proposition~\ref{LemmaplantedCluster}.
The following statement relates the planted model to the random
hypergraph $H_k(n,m)$.
{A similar statement has been obtained independently by Achlioptas and
Theodoropoulos~\cite{AchlioptasTheodoropoulos}.}

%
\begin{lemma}\label{LemmaantiPlanting}
Let $d>0$ and $\beta\ge0$. Assume that there is a sequence
$({{\mathcal E}}
_n)_{n\geq1}$ of events such that
$\limsup_{n\rightarrow\infty}\pr\brk{\G\in{{\mathcal E}}_n}^{1/n}<1$.
Then $\Erw[\Zkc(H_k(n,\break m))\vecone_{{{\mathcal E}}_n}]\le
\exp
(-\Omega
(n))\Erw[\Zkc(H_k(n,m))]$.
\end{lemma}

\begin{pf}
Fix $\alpha>0$ such that $\limsup_{n\rightarrow\infty}\pr\brk{\G
\in{{\mathcal E}}
_n}^{1/n}\le\exp(-\alpha)$. To shorten the notation, we write
$H_{n,m}$ for $H_k(n,m)$. For any $\eps>0$, we have the decomposition
%
%
\begin{eqnarray}\label{eqnewlabelsec50}
&& \Erw\bigl[\Zkc(H_{n,m})\vecone_{{{\mathcal E}}_n}\bigr]\nonumber
\\
&&\qquad = \sum
_{\sigma:\brk n\rightarrow\{{\pm1} \}}\Erw\bigl[\exp\bigl( -
\beta E_{H_{n,m}}(\sigma)\bigr)\vecone_{{{\mathcal E}}_n}\bigr]
\\
&&\qquad \leq\sum_{\sigma\in B_\eps}\Erw\bigl[\exp\bigl( -\beta
E_{H_{n,m}}(\sigma)\bigr)\vecone_{{{\mathcal E}}_n}\bigr] +  \sum_{\sigma\notin B_\eps}\Erw
\bigl[\exp\bigl( -\beta E_{H_{n,m}}(\sigma)\bigr)\bigr].\nonumber
\end{eqnarray}
To bound the first summand in (\ref{eqnewlabelsec50}), we let
$m_0=\frac{2^{1-k}\exp(-\beta)}{1-2^{1-k}(1-\exp(-\beta))}m$ and
define the set $M_\eps= \{{\mu\in\brk m:\mid\mu-m_0\mid<\eps
n} \}$.
Now, for any $\mu\in[m]$ we have
\begin{eqnarray*}
&&\sum_{\sigma\in B_\eps}\pr\bigl[\bigl\{E_{H_{n,m}}(
\sigma)=\mu\bigr\} \cap\{ H_{n,m} \in{{\mathcal E}}_n\}
\bigr]
\\
&&\qquad =\sum_{\sigma\in B_\eps}\pr\bigl[H_{n,m} \in{{
\mathcal E}}_n\mid E_{H_{n,m}}(\sigma)=\mu\bigr]\pr
\brkk{E_{H_{n,m}}(\sigma)=\mu}.
\end{eqnarray*}
Under the conditions $e(\G)=m$ and $E_{H_{n,m}}(\sigma)=E_{\G
}(\sigma)$ for $\sigma:\brk n\rightarrow\{{\pm1} \}$,
the two random
hypergraphs $H_{n,m}$ and $\G$ are identically distributed. Therefore,
\begin{eqnarray*}
&& \pr\bigl[H_{n,m} \in{{\mathcal E}}_n\bigl\vert
E_{H_{n,m}}(\sigma)=\mu\bigr] \nonumber
\\
&&\qquad =\pr\bigl[\G\in{{\mathcal E}}_n
\bigr\vert E_{\G}(\sigma)=\mu, e(\G)=m\bigr]
\leq\frac{\pr[\G\in{{\mathcal E}}_n]}{\pr\brkk{E_{\G}(\sigma
)=\mu,e(\G)=m}}.
\end{eqnarray*}
By standard concentration results, there is $\eps>0$ such that
\[
\pr\brkk{E_{\G}(\sigma)=\mu, e(\G)=m}\geq\exp\biggl(-
\frac
{\alpha}2 n \biggr)\qquad\mbox{for any }\sigma\in B_\eps,\mu
\in M_\eps.
\]
Hence, for any $\mu\in M_\eps$:
\begin{eqnarray*}
&&\sum_{\sigma\in B_\eps}\pr\bigl[\bigl\{E_{H_{n,m}}(
\sigma)=\mu\bigr\}\cap\{ H_{n,m} \in{{\mathcal E}}_n\}\bigr]
\\
&&\qquad \leq \exp\biggl(\frac{\alpha}2 n \biggr)\sum_{\sigma\in B_\eps
}
\pr[\G\in{{\mathcal E}}_n]\pr\brkk{E_{H_{n,m}}(\sigma)=\mu}
\end{eqnarray*}
and, therefore, letting $A=2^n(1-2^{1-k}(1-\exp(-\beta)))^m$, we get
\begin{eqnarray}\label{eqnewlabelsec53}
&& \sum_{\mu\in M_\eps}\sum_{\sigma\in B_\eps}
\Erw\bigl[\exp\bigl(-\beta E_{H_{n,m}}(\sigma)\bigr)\vecone
_{{{\mathcal E}}_n}
\bigr]
\nonumber
\\
&&\qquad =\sum_{\mu\in M_\eps}\sum_{\sigma\in B_\eps}
\exp(-\beta\mu)\pr\bigl[\bigl\{E_{H_{n,m}}(\sigma)=\mu\bigr\}
\cap
\{H_{n,m} \in{{\mathcal E}}_n\} \bigr]
\nonumber\\[-8pt]\\[-8pt]\nonumber
&&\qquad \leq\exp\biggl(-\frac{\alpha}2n \biggr)\sum_{\mu\in M_\eps}
\sum_{\sigma\in B_\eps}\exp(-\beta\mu)\pr\brkk{E_{H_{n,m}}(
\sigma)=\mu}
\nonumber
\\
&&\qquad \leq A\exp\biggl(-\frac{\alpha}2n \biggr).\nonumber
\end{eqnarray}
Furthermore, Lemma~\ref{LemavgEngergy} shows that there is $\delta
>0$ such that
%
\begin{equation}
\label{eqnewlabelsec52} \sum_{\mu\notin M_\eps
}\sum
_{\sigma\in B_\eps}\exp(-\beta\mu)\pr\bigl[E_{H_{n,m}}(\sigma
)=\mu
\bigr]\leq A\exp(-\delta n).
\end{equation}
To bound the second summand in (\ref{eqnewlabelsec50}), we get
from Lemma~\ref{Lemimbalanced} that there is $\delta'>0$ such that
%
\begin{equation}
\label{eqnewlabelsec51} \sum_{\sigma\notin
B_\eps}\Erw\bigl[\exp\bigl(-\beta
E_{H_{n,m}}(\sigma)\bigr)\bigr]\leq A\exp\bigl(-\delta' n
\bigr).
\end{equation}
Combining the estimates (\ref{eqnewlabelsec53}), (\ref
{eqnewlabelsec52}) and (\ref{eqnewlabelsec51}) in the
decomposition (\ref{eqnewlabelsec50})
yields
\[
\Erw\bigl[\Zkc(H_{n,m})\vecone_{{{\mathcal E}}_n}\bigr]\le A\exp
\bigl( -
\max\bigl(\alpha/2,\delta,\delta'\bigr)n\bigr).
\]
The assertion follows with Lemmas~\ref{LemmafirstMoment} and~\ref{Lemrelatepartition}.
\end{pf}

%
\begin{corollary} \label{CorollaryantiPlanting}
Let $d>0$ and $\beta\ge0$. Assume that there exists a sequence
$({{\mathcal E}}
_n)_{n\geq1}$ of events such that
\[
\lim_{n\rightarrow\infty}\pr\brkk{H_k(n,m)\in{{\mathcal
E}}_n}=1\qquad\mbox{while } \limsup_{n\rightarrow\infty}\pr\brk
{\G
\in{{\mathcal E}}_n}^{1/n}<1.
\]
Then $\Phi_{d,k}(\beta)<\ln2+\frac{d}k\ln(1-2^{1-k}(1-\exp(-\beta)))$.
\end{corollary}
\begin{pf}
Since $Z_\beta(H_k(n,m))^{1/n} \leq2$ and $\mathbb{P} [
H_k(n,m) \in{{\mathcal E}}_n ] = 1-o(1)$, Jensen's inequality yields
\begin{eqnarray*}
\Erw\bigl[ Z_\beta\bigl(H_k(n,m)\bigr)^{1/n}
\bigr]& =& \Erw\bigl[ Z_\beta\bigl(H_k(n,m)
\bigr)^{1/n} \vecone_{{{\mathcal E}}_n} \bigr] +o(1)
\\
&\leq&\Erw\bigl[ Z_\beta\bigl(H_k(n,m)\bigr)
\vecone_{{{\mathcal E}}_n} \bigr]^{1/n} +o(1).
\end{eqnarray*}

Hence, under the assumptions of the corollary~we obtain with Jensen's
inequality and Lemma~\ref{LemmaantiPlanting}
\begin{eqnarray*}
\Phi_{d,k}(\beta)&\le&\limsup_{n \to\infty}\ln\Erw\bigl[
Z_\beta\bigl(H_k(n,m)\bigr)^{1/n} \bigr]
\\
& \le& \exp\bigl(- \Omega(1)\bigr) \limsup_{n \to\infty}\Erw
\bigl[
Z_\beta\bigl(H_k(n,m)\bigr) \bigr]^{1/n}.
\end{eqnarray*}
The result then follows from Lemmas~\ref{LemmafirstMoment} and~\ref{Lemrelatepartition}.
\end{pf}

\subsection{An unlikely event}\label{Secpartition}
As a next step, we establish the following.

%
\begin{lemma}\label{Lemmapartition}
Assume that~(\ref{eqLemmaplantedCluster3}) holds for some $\beta\ge
k\ln2-\ln k$.
Then there exists $z>0$ such that
\[
\lim_{n\rightarrow\infty}\pr\brkkk{\frac{1}n\ln Z_{\beta
}
\bigl(H_k(n,m)\bigr)\leq z}=1,\qquad\limsup_{n\rightarrow\infty}
\pr\brkkk{\frac{1}n\ln Z_{\beta}(\G)\leq z}^{1/n}<1.
\]
\end{lemma}

The proof of Lemma~\ref{Lemmapartition}, to which we dedicate the
rest of this subsection, is an extension of the argument from~\cite
{bapst}, Section~6, to the case of finite $\beta$.
We need the following concentration result.

%
\begin{lemma}\label{CorZAzuma}
For any fixed $d>0$, $\beta\ge0$, $\alpha>0$ there are $\delta>0$,
$\delta'>0$ such that the following is true.
Suppose that $(\sigma_n)_{n\geq1}$ is a sequence of maps $\brk
n\rightarrow
\{{\pm1} \}$.
Then for all large enough $n$,
\[
\pr\brkk{\bigl\vert\ln\bigl(Z_{\beta}(\G)\bigr)-\Erw\bigl[\ln
Z_{\beta}(\G)\bigr\vert\bolds{\sigma}=\sigma_n\bigr]\vert>
\alpha n\vert\bolds{\sigma}=\sigma_n}\leq\exp(-\delta n)
\]
and
\[
\pr\brkk{\bigl\vert\ln\bigl({{\mathcal C}}_{\beta}(\G,\bolds
{\sigma})\bigr)-
\Erw\bigl[\ln{{\mathcal C}}_{\beta}(\G,\bolds{\sigma})\bigr
\vert\bolds{\sigma}=
\sigma_n\bigr]\vert>\alpha n\vert\bolds{\sigma}=\sigma_n}
\leq\exp\bigl(-\delta' n\bigr).
\]
\end{lemma}
\begin{pf}
This is immediate from the Lipschitz property and McDiarmid's
inequality~\cite{McDiarmid}, Theorem~3.8.
\end{pf}

We further need several statements about quantities in the planted
model conditioned on $\bolds{\sigma}$ being some fixed (balanced) coloring.

%
\begin{lemma}\label{LemmapickAndChoose}
Assume that~(\ref{eqLemmaplantedCluster3}) is true for some $\beta
\ge k\ln2-\ln k$.
Then there exist a fixed number $\eps>0$ and a sequence $\sigma_n$ of
balanced maps $\brk n\rightarrow\{{\pm1} \}$ such that
\begin{eqnarray*}
&&\lim_{n\rightarrow\infty}\pr\biggl[\frac{1}n\ln{{\mathcal
C}}_\beta(\G, \bolds{\sigma}) > \ln2 +
\frac{d}k\ln\bigl(1 - 2^{1-k}
\bigl(1 - \exp(-\beta)\bigr)\bigr)+\eps\mid\bolds{\sigma}=
\sigma_n \biggr] = 1.
\end{eqnarray*}
\end{lemma}

\begin{pf}
By Stirling's formula, there is an $n$-independent number $\delta>0$
such that for sufficiently large $n$ we have
%
%
\begin{equation}
\label{eqLemmapickAndChoose} \pr\brk{\bolds{\sigma}\in
\Bal}\geq\delta.
\end{equation}
Let $A=\ln2+\frac{d}k\ln(1-2^{1-k}(1-\exp(-\beta)))$. Using (\ref
{eqLemmaplantedCluster3}), we know that there is $\eps>0$ such that
$\liminf_{n\rightarrow\infty}\pr\brk{\frac{1}n\ln{{\mathcal
C}}_\beta(\G, \bolds{\sigma}
)>A+3\eps}\ge0.9$. With the concentration bound from Lemma~\ref
{LemmaMcDiarmid}, we get
\[
\lim_{n\rightarrow\infty}\pr\brkkk{\frac{1}n\ln{{\mathcal
C}}_\beta(\G, \bolds{\sigma})>A+2\eps}=1.
\]
Thus, with $p_n=\liminf_{n\rightarrow\infty}\max_{\sigma_n \in
\Bal}\pr
\brk{\frac{1}n\ln{{\mathcal C}}_\beta(\G, \bolds{\sigma
})>A+2\eps\mid\bolds{\sigma}
=\sigma
_n}$ and (\ref{eqLemmapickAndChoose}) we get
\begin{eqnarray}
1&\le&\liminf_{n\rightarrow\infty} \biggl(\sum_{\sigma_n\in\Bal
} \pr\brkkk{\frac{1}n\ln{{\mathcal C}}_\beta(\G, \bolds{\sigma
})>A+2\eps
\mid\bolds{\sigma}=\sigma_n}\pr[\bolds{\sigma}=\sigma_n]\nonumber
\\
&&{}+\sum
_{\sigma_n \notin\Bal} \pr[\bolds{\sigma}=\sigma_n] \biggr)
\nonumber\\[-8pt]\\[-8pt]\nonumber
&\le&\liminf_{n\rightarrow\infty}p_n \pr[\bolds{\sigma}\in\Bal
]+\pr[
\bolds{\sigma}\notin\Bal]
\nonumber
\\
& \le&\liminf_{n\rightarrow\infty}p_n +1-\delta,
\nonumber
\end{eqnarray}
implying that $\liminf_{n\rightarrow\infty}p_n\ge\delta$. Thus, the
concentration bound from Lem\-ma~\ref{CorZAzuma} yields
\[
\lim_{n\rightarrow\infty}\max_{\sigma_n \in\Bal}\pr\brkkk{\frac{1}n\ln{{\mathcal C}}_\beta(\G, \bolds{\sigma})>A+\eps\mid
\bolds{\sigma}=
\sigma_n}=1
\]
completing the proof.
\end{pf}

%
\begin{lemma}\label{LemmanearlyBalanced}
For any $\eta>0$, there is $\delta>0$ such that
\[
\limsup_{n\rightarrow\infty}\frac{1}n\ln\pr\brkk{\bigl\vert
\bigl
\vert\bolds{\sigma}^{-1}(1)\bigr\vert-n/2\bigr\vert>\eta n}
\leq-\delta.
\]
\end{lemma}
\begin{pf}
This is immediate from the Chernoff bound.
\end{pf}

For a set $S \subset V$ let $\Vol(S\mid H)$ be the sum of the degrees of
the vertices in $S$ in the hypergraph $H$.

%
\begin{lemma}\label{LemmaVol}
For any $\gamma>0$, there is $\alpha>0$ such that for any set
$S\subset\brk n$
of size $\vert S\vert\leq\alpha n$ and any map $\sigma:\brk
n\rightarrow
\{{\pm1} \}$
we have
$\limsup\frac{1}n\ln\pr[ \Vol(S\mid\G)\ge\gamma n\mid
\bolds{\sigma}=\sigma
]\leq-\alpha$.
\end{lemma}
\begin{pf}
Let $(X_v)_{v\in\brk n}$ be a family of independent random variables
with distribution $\Bin(\bink{n-1}{k-1},2p )$.
Then for any $\sigma$ and any $S\subset\brk n$ the volume $\Vol
(S\mid\G
)$ is stochastically dominated by $X_S=2k\sum_{v\in S}X_v$.
Furthermore, $\Erw[X_S]= 4dk\vert S\vert$.
Thus, for any $\gamma>0$ we can choose an $n$-independent $\alpha>0$
such that for any $S\subset\brk n$ of size $\vert S\vert\leq\alpha
n$ we have
$\Erw[X_S]\leq\gamma n/2$.
In fact, the Chernoff bound shows that by picking $\alpha>0$
sufficiently small, we can ensure that
$\pr\brk{\Vol(S\mid\G)\geq\gamma n\mid\bolds{\sigma}=\sigma
}\leq\pr\brk
{X_S\geq\gamma n}\leq\exp(-\alpha n)$,
as desired.
\end{pf}

%
\begin{lemma}\label{LemmasigmaRandom}
Let $d>0$ and $\beta\ge0$. Assume that there exist numbers $z>0$,
$\eps>0$ and a sequence $(\sigma_n)_{n\geq1}$ of balanced maps $\brk
n\rightarrow\{{\pm1} \}$ such that
\[
\lim_{n\rightarrow\infty}\frac{1}n\Erw\brkk{\ln Z_{\beta}(\G)
\mid\bolds{\sigma}=\sigma_n}>z+\eps.
\]
Then
$\limsup_{n\rightarrow\infty}\pr\brk{\frac{1}n\ln Z_{\beta}(\G
)\leq z}^{1/n}<1$.
\end{lemma}
\begin{pf}
Suppose that $n$ is large enough so that $\frac{1}n\Erw\brk{\ln
Z_{\beta}(\G)\mid\bolds{\sigma}=\sigma_n}>z+\eps/2$.
Set $n_i=\vert\sigma_n^{-1}(i)\vert$ and let $T$ be the set of all
$\tau:\brk
n\rightarrow\{{\pm1} \}$ such that $\vert\tau^{-1}(i)\vert=n_i$
for $i=\pm1$.
As $Z_{\beta}$ is invariant under permutations of the vertices, we have
%
%
\begin{eqnarray}\label{eqLemmasigmaRandom1}
\frac{1}n\Erw\brkk{\ln Z_{\beta}(\G
)\mid\bolds{\sigma}=
\tau}=\frac{1}n\Erw\brkk{\ln Z_{\beta}(\G)\mid\bolds{\sigma}=
\sigma_n}>z+\eps/2
\nonumber\\[-8pt]\\[-8pt]
\eqntext{\mbox{for any }\tau\in T.}
\end{eqnarray}
Let $\gamma=\eps/(4\beta)>0$.
By Lemma~\ref{LemmaVol}, there exists $\alpha>0$ such that for large
enough $n$ for any set $S\subset V$ of size $\vert S\vert\leq\alpha
n$ and
any $\sigma:\brk n\rightarrow\{{\pm1} \}$ we have
%
%
\begin{equation}
\label{eqLemmasigmaRandom4} \pr\brkkk{\Vol(S\mid\G)<\frac{\gamma
n}2\Big| \bolds{\sigma}=\sigma}\geq1-
\exp(-\alpha n).
\end{equation}
Fix such an $\alpha>0$, and pick and fix a small $0<\eta<\alpha/3$.
By Lemma~\ref{LemmanearlyBalanced}, there exists an ($n$-independent)
number $\delta=\delta(\beta,\eps,\eta)>0$ such that
%
%
\begin{equation}
\label{eqLemmasigmaRandom3} \pr\brk{\bolds{\sigma}\in B_\eta}\geq
1-\exp(-\delta n).
\end{equation}
Because $\sigma_n$ is balanced, we have $\vert n_i-n/2\vert\leq
\sqrt n$ for
$i=\pm1$. Therefore, if $\bolds{\sigma}\in B_\eta$, then it is
possible to
obtain from $\bolds{\sigma}$ a map $\tau_{\bolds
{\sigma}}\in T$ by
changing the colors of at most $2\eta n$ vertices.
Hence, if $\bolds{\sigma}\in B_\eta$ we let $\G_{\tau_{\bolds
{\sigma}}}$ be the
random hypergraph with planted coloring $\tau_{\bolds{\sigma}}$.
Further, let $\G_{\bolds{\sigma}}$ be the hypergraph obtained by removing
from $\G_{\tau_{\bolds{\sigma}}}$ each edge
that is monochromatic under $\bolds{\sigma}$ but not under $\tau
_{\bolds{\sigma}}$
with probability $1-\exp(-\beta)$ independently and inserting each
edge that is monochromatic under $\tau_{\bolds{\sigma}}$ but not under
$\bolds{\sigma}$ with probability $(1-\exp(-\beta))p_2$ independently.
Then $\G_{\bolds{\sigma}}=\G$ in distribution.

Let $S_{\bolds{\sigma}}$ be the set of vertices $v$ with
$\bolds{\sigma}(v)\neq\tau_{\bolds{\sigma}}(v)$. Our
choice of $\eta$
ensures that $\vert S_{\bolds{\sigma}}\vert<\alpha n$. Let $\Delta
$ be the number of
edges present in $\G_{\tau_{\bolds{\sigma}}}$ but not in $\G
_{\bolds{\sigma}}$ or
vice versa.
Then $\Delta\leq\Vol(S_{\bolds{\sigma}}\vert\G_{\tau
_{\bolds{\sigma}}})+\Vol
(S_{\bolds{\sigma}}\vert\G_{\bolds{\sigma}})$.
Hence, with (\ref{eqLemmasigmaRandom4}) there exists a constant $c>0$
such that
%
%
\begin{equation}
\label{eqLemmasigmaRandom4a} \pr\brk{\Delta\leq\gamma n\mid\bolds
{\sigma}\in B_\eta}\geq1-c
\exp(-\alpha n).
\end{equation}
Using (\ref{eqLemmasigmaRandom3}), (\ref{eqLemmasigmaRandom4a}) and
the fact that removing a single edge can reduce $\frac{1}n \ln
Z_{\beta
}$ by at most $\beta/n$, we obtain
\begin{eqnarray}\label{eqLemmasigmaRandom5}
\pr\biggl[\frac{1}n \ln Z_{\beta}(\G)\leq z \biggr]&=&\pr\biggl[
\frac{1}n \ln Z_{\beta}(\G_{\bolds{\sigma}})\leq z \biggr]
\nonumber
\\
&\leq&\exp(-\delta n)+\pr\biggl[\frac{1}n \ln Z_{\beta}(
\G_{\bolds{\sigma}
})\leq z\Big| \bolds{\sigma}\in B_\eta\biggr]
\nonumber
\\
&\leq&\exp(-\delta n)+c\exp(-\alpha n)
\nonumber\\[-8pt]\\[-8pt]\nonumber
&&{} +\pr\biggl[\frac{1}n \ln Z_{\beta}(\G_{\bolds{\sigma}})\leq z
\Big|  \bolds{\sigma}\in B_\eta,\Delta\leq\gamma n \biggr]
\nonumber
\\
&\leq&\exp(-\delta n)+c\exp(-\alpha n)
\nonumber
\\
&&{} +\pr\biggl[\frac{1}n \ln Z_{\beta}(\G_{\tau_{\bolds{\sigma
}}})-\gamma
\beta\leq z\Big| \bolds{\sigma}\in B_\eta,\Delta\leq\gamma n
\biggr].\nonumber
\end{eqnarray}
By the choice of $\gamma$, (\ref{eqLemmasigmaRandom3}), (\ref
{eqLemmasigmaRandom4a}) and (\ref{eqLemmasigmaRandom1}), we have
\begin{eqnarray*}
&&\pr\biggl[\frac{1}n \ln Z_{\beta}(\G_{\tau_{\bolds{\sigma
}}})-\gamma
\beta\leq z\Big|  \bolds{\sigma}\in B_\eta,\Delta\leq\gamma n
\biggr]
\\
&&\qquad  \le2\pr\biggl[\frac{1}n \ln Z_{\beta}(\G_{\tau_{\bolds{\sigma
}}})\leq
z+\frac{\eps}{4}\Big|  \bolds{\sigma}\in B_\eta\biggr]
\\
&&\qquad \le3\pr\biggl[\frac{1}n \ln Z_{\beta}(\G)\leq z+
\frac{\eps
}{4}\Big|  \bolds{\sigma}=\sigma_n \biggr]
\\
&&\qquad  \le3\pr\biggl[\frac{1}n \ln Z_{\beta}(\G)\leq\frac{1}n
\Erw\bigl[\ln Z_{\beta}(\G)\mid \bolds{\sigma}=\sigma_n
\bigr]-\frac{\eps}{4}\Big|  \bolds{\sigma}=\sigma_n \biggr].
\end{eqnarray*}
The assertion follows by combining this with (\ref
{eqLemmasigmaRandom5}) and Lemma~\ref{CorZAzuma}.
\end{pf}

\begin{pf*}{Proof of Lemma~\ref{Lemmapartition}}
Lemma~\ref{LemmapickAndChoose} shows that there exist $\eps>0$ and
balanced maps $\sigma_n:\brk n\rightarrow\{{\pm1} \}$
such that
%
%
\begin{eqnarray}\label{eqLemmapartition1}
&& \lim_{n\rightarrow\infty}\pr\biggl[\frac{1}n\ln{{\mathcal
C}}_\beta(\G,\bolds{\sigma}) \geq\ln2
+ \frac{d}k\ln\bigl(1 - 2^{1-k}
\bigl(1 - \exp(-\beta)\bigr)\bigr)+\eps\Big|  \bolds{\sigma}=
\sigma_n \biggr] \hspace*{-25pt}
\nonumber\\[-8pt]\\[-8pt]\nonumber
&&\qquad = 1.
\end{eqnarray}
Clearly, (\ref{eqLemmapartition1}) implies that
%
%
\begin{eqnarray}\label{eqLemmapartition2}
\qquad && \lim_{n\rightarrow\infty}\pr\biggl
[\frac{1}n\ln
Z_\beta(\G) \geq\ln2 +
\frac{d}k\ln\bigl(1 - 2^{1-k}\bigl(1
- \exp(-\beta)\bigr)\bigr)+\eps\Big|  \bolds{\sigma}=
\sigma_n \biggr]
\nonumber\\[-8pt]\\[-8pt]\nonumber
&&\qquad  = 1.
\end{eqnarray}
Hence, with $z=\ln2+\frac{d}k\ln(1-2^{1-k}(1-\exp(-\beta)))+\eps
/2$, Lemma~\ref{LemmasigmaRandom} and~(\ref{eqLemmapartition2}) yield
%
%
\begin{equation}
\label{eqLemmapartition2b} \limsup_{n\rightarrow\infty}\pr\brkkk
{\frac{1}n\ln
Z_{\beta}(\G)\leq z}^{1/n}<1.
\end{equation}
By comparison, Lemma~\ref{LemmafirstMoment} and Lemma~\ref
{Lemrelatepartition} imply
%
%
\begin{equation}
\label{eqLemmapartition3} \lim_{n\rightarrow\infty}\pr\brkkk{\frac
{1}n\ln
Z_{\beta
}\bigl(H_k(n,m)\bigr)\leq z}=1.
\end{equation}
Thus, the assertion follows from (\ref{eqLemmapartition2b})
and~(\ref{eqLemmapartition3}).
\end{pf*}

\subsection{Tame colorings}

To facilitate the proof of the first part of Proposition~\ref
{LemmaplantedCluster}, we introduce a random variable that explicitly
controls the ``cluster size'' ${{\mathcal C}}_\beta(H_k(n,m),\sigma)$.
The idea of explicitly controlling the cluster size was introduced
in~\cite{Lenka} in the ``zero temperature'' case, and here
we generalise it to the case of finite $\beta$.
More precisely, we call $\sigma:\brk n\rightarrow\{{\pm
1} \}$ \textit{tame} in
$H$ if $\sigma$ is balanced and if
${{\mathcal C}}_\beta(H,\sigma)\leq\Erw[Z_\beta(H)]$.
Now, let
\[
\Zkg\bigl(H_k(n,m)\bigr)=\sum_{\sigma:\brk n\rightarrow\{{-1,1} \}
}\exp
\bigl(-\beta E_{H_k(n,m)
}(\sigma)\bigr)\cdot\vecone_{\sigma \mbox{\,\fontas is tame}}.
\]

%
\begin{lemma}\label{Lemmasmm}
Assume that $0\leq d/k\leq2^{k-1}\ln2+O_k(1)$ is such that\break 
$\liminf_{n\rightarrow\infty}\frac{\Erw[\Zkg(H_k(n,m))]}{\Erw
[Z_\beta(H_k(n,m))]}>0$.
Then
\[
\liminf_{n\rightarrow\infty}\frac{\Erw[\Zkg(H_k(n,m))]^2}{\Erw
[\Zkg
(H_k(n,m))^2]}>0.
\]
\end{lemma}
\begin{pf}
The proof is based on a second moment argument.
Mimicking the notation of Section~\ref{Secsecond}, we let
\begin{eqnarray*}
&& \Zkg(\alpha)
\\
&&\qquad =\sum_{\sigma,\tau:\scal{\sigma}{\tau}=\alpha n}\exp\bigl
(-\beta\bigl(E_{H_k(n,m)}(
\sigma)+E_{H_k(n,m)}(\tau)\bigr)\bigr)\cdot\vecone_{\sigma\mbox
{\,\fontas is tame}} \cdot
\vecone_{\tau\mbox{\,\fontas is tame}}.
\end{eqnarray*}
Then it is clear that
\[
\Erw\bigl[\Zkg\bigl(H_k(n,m)\bigr)^2\bigr]=
\sum_{\nu=0}^n\Erw\brkk{\Zkg(2\nu/n-1)}.
\]
Furthermore, we have $\Zkg(\alpha)\leq\Zkc(\alpha)$ for any
$\alpha$.
We define $I=[-1+2^{-3k/4},1-2^{-3k/4}]$.
Lemma~\ref{Lemmasmm1} and Lemma~\ref{LemmaLaplace} yield
%
%
\begin{equation}
\label{eqLemmasmm3} \sum_{\alpha\in I}\Erw\brkk{\Zkc(\alpha)}= O
\bigl(\Erw\bigl[\Zkc\bigl(H_k(n,m)\bigr)\bigr]^2\bigr).
\end{equation}
By the definition of ``tame'' we have
%
%
\begin{eqnarray}\label{eqLemmasmm4}
&&\sum_{\alpha>1-2^{-3k/4}}\Erw\brkk{\Zkg(\alpha)}\nonumber
\\
&&\qquad  \leq\Erw\brkkkkp{\sum_{\sigma}\exp\bigl(-\beta
E_{H_k(n,m)}(\sigma)\bigr)\cdot\vecone_{\sigma\mbox{\,\fontas
is tame}} \cdot{{\mathcal
C}}_\beta\bigl(H_k(n,m),\sigma\bigr)}
\nonumber\\[-8pt]\\[-8pt]\nonumber
&&\qquad  \le\Erw\brkkkkp{\sum_{\sigma}\exp\bigl(-\beta
E_{H_k(n,m)}(\sigma)\bigr) \cdot\Erw\brkk{\Zkg\bigl(H_k(n,m)
\bigr)}}
\\
&&\qquad  = O\bigl(\Erw\brkk{\Zkg\bigl(H_k(n,m)\bigr)}^2\bigr).
\nonumber
\end{eqnarray}
Moreover, $\sum_{\alpha<-1+2^{-3k/4}}\Erw\brk{\Zkg(\alpha)}=\sum
_{\alpha>1-2^{-3k/4}}\Erw\brk {\Zkg(\alpha)}$ by symmetry.
Hence, $\Erw[\Zkg(H_k(n,m))^2]= O(\Erw[\Zkc(H_k(n,m))]^2)$ by equations
(\ref{eqLemmasmm3}) and (\ref{eqLemmasmm4}).

Finally, the assertion follows from our assumption that
$\Erw[\Zkg(H_k(n,m))]=\Omega(\Erw[\Zkc(H_k(n,\break m))])$.
\end{pf}

%
\begin{lemma}\label{LemmatamePlanting}
Let $d>0$ and $\beta\ge0$ and assume that we have
\[
\limsup
_{n\rightarrow\infty}\pr\brk{\bolds{\sigma}\mbox{ is not tame in
} \G}^{1/n}<1.
\]
Then there is $c>0$ such that $\Erw[\Zkg(H_k(n,m))]\ge\Erw[\Zkc
(H_k(n,m))]/c$.
\end{lemma}

\begin{pf}
The proof is very similar to the proof of~Lemma~\ref
{LemmaantiPlanting}. We fix an $\alpha>0$ such that $\limsup
_{n\rightarrow
\infty}\pr\brk{\bolds{\sigma}\mbox{  is not tame in } \G
}^{1/n}\le\exp
(-\alpha)<1$.
For any $\eps>0$, we have
\begin{eqnarray*}
&& \Erw\bigl[\Zkc\bigl(H_k(n,m)\bigr)-\Zkg\bigl(H_k(n,m)
\bigr)\bigr]
\\
&&\qquad  =\sum_{\sigma:\brk n\rightarrow\{{\pm1} \}}\Erw\bigl[\exp\bigl
(-\beta
E_{H_k(n,m)
}(\sigma)\bigr)\vecone_{\sigma\mbox{\,\fontas is not tame in }
H_k(n,m)}\bigr]
\\
&&\qquad  \leq\sum_{\sigma\in B_\eps}\Erw\bigl[\exp\bigl(-\beta
E_{H_k(n,m)}(\sigma)\bigr)\vecone_{\sigma\mbox{\,\fontas is not tame in }
H_k(n,m)}\bigr]
\\
&&\quad\qquad  {} +\sum_{\sigma\notin B_\eps}\Erw\bigl[\exp\bigl(-\beta
E_{H_k(n,m)}(\sigma)\bigr)\bigr].
\end{eqnarray*}
With $m_0$ and $M_\eps$ as in the proof of~Lemma~\ref
{LemmaantiPlanting} and $\mathcal A(\sigma,\mu)$ the event $\{E_{\G
}(\sigma)=\mu, e(\G)=m, \vert\sigma^{-1}(1)\vert=\vert\bolds
{\sigma}^{-1}(1)\vert\}$, we
fix an $\eps>0$ such that\break $\pr[\mathcal A(\sigma,\mu)]>\exp
(-\frac{\alpha}2n )$ for all $\sigma\in B_\eps,\mu\in M_\eps
$. Then for any $\mu\in M_\eps$:
\begin{eqnarray*}
&&\sum_{\sigma\in B_\eps}\pr\bigl[\bigl\{E_{H_k(n,m)}(
\sigma)=\mu\bigr\} \cap\bigl\{ \sigma\mbox{ is not tame in }H_k(n,m)
\bigr\}\bigr]
\\
&&\qquad  =\sum_{\sigma\in B_\eps}\pr\bigl[\sigma\mbox{ is not tame in
}H_k(n,m) \mid E_{H_k(n,m)}(\sigma)=\mu\bigr]\pr
\brkk{E_{H_k(n,m)}(\sigma)=\mu}
\\
&&\qquad  =\sum_{\sigma\in B_\eps}\pr\bigl[\bolds{\sigma}\mbox{ is not tame in }
\G\mid\mathcal A(\sigma,\mu)\bigr]\pr\brkk{E_{H_k(n,m)}(\sigma
)=\mu}
\\
&&\qquad  \leq\sum_{\sigma\in B_\eps}\frac{\pr[\bolds{\sigma}\mbox{ is not tame in }\G]}{\pr(\mathcal A(\sigma,\mu))}\pr
\brkk{E_{H_k(n,m)
}(\sigma)=\mu}
\\
&&\qquad  \leq\exp\biggl(-\frac{\alpha}2n \biggr) \sum
_{\sigma\in B_\eps
} \pr\brkk{E_{H_k(n,m)}(\sigma)=\mu}.
\end{eqnarray*}
Letting $A=2^n(1-2^{1-k}(1-\exp(-\beta)))^m$, we get
%
%
\begin{eqnarray}
\label{eqtame1} &&\sum_{\mu\in M_\eps}\sum
_{\sigma\in B_\eps}\Erw\bigl[\exp\bigl(-\beta E_{H_k(n,m)}(\sigma)
\bigr)\vecone_{\sigma\mbox{\,\fontas is not tame in }
H_k(n,m)}\bigr]\nonumber
\\
&&\qquad =\sum_{\mu\in M_\eps}\sum
_{\sigma\in B_\eps}\exp(-\beta\mu)\pr\bigl[\bigl\{E_{H_k(n,m)}(
\sigma)=\mu\bigr\}
\\
\nonumber
&&\quad\qquad {} \cap\bigl\{\sigma\mbox{ is not tame in }H_k(n,m)\bigr\}
\bigr]\leq A\exp\biggl(-\frac{\alpha}2n \biggr).
\end{eqnarray}
Furthermore Lemma~\ref{LemavgEngergy} shows that there is $\delta>0$
such that
%
%
\begin{eqnarray}
&&\label{eqtame2} \sum_{\mu\notin M_\eps}\sum
_{\sigma\in B_\eps}\exp(-\beta\mu)\pr\bigl[E_{H_k(n,m)}(\sigma
)=\mu
\bigr]\leq A\exp(-\delta n)
\end{eqnarray}
and we get from Lemma~\ref{Lemimbalanced} that there is $\delta'>0$
such that
%
%
\begin{eqnarray}
\label{eqtame3}
&&\sum_{\sigma\notin B_\eps}\Erw\bigl[\exp\bigl(-
\beta E_{H_k(n,m)}(\sigma)\bigr)\bigr]\leq A\exp\bigl(-\delta'
n\bigr).
\end{eqnarray}
Combining the estimates (\ref{eqtame1}), (\ref{eqtame2}) and
(\ref{eqtame3}) and using Lemmas~\ref{LemmafirstMoment} and~\ref{Lemrelatepartition}
yields
\begin{eqnarray*}
\Erw\bigl[\Zkc\bigl(H_k(n,m)\bigr)-\Zkg\bigl(H_k(n,m)
\bigr)\bigr]&\le& A\exp\bigl( - \max\bigl(\alpha/2,\delta,\delta'
\bigr)n\bigr)
\\
& \le&\exp\bigl(-\Omega(n)\bigr)\Erw\bigl[\Zkc\bigl
(H_k(n,m)\bigr)
\bigr],
\end{eqnarray*}
which proves the assertion.
\end{pf}

%
\begin{corollary}\label{LemmaplantedClustereasy}
Assume that $d/k=2^{k-1}\ln2+O_k(1)$ and that $\beta_0 \ge k\ln2-\ln
k$ is such that~(\ref{eqLemmaplantedCluster2}) holds for all $k\ln
2-\ln k\leq\beta\leq\beta_0$.
Then $\beta_{\mathrm{crit}}(d,k)\geq\beta_0$.
\end{corollary}

The proof of this corollary~extends a ``zero temperature'' argument
from~\cite{bapst}, Section~5, to the case of $\beta\in[0,\infty)$.

\begin{pf*}{Proof of Corollary \ref{LemmaplantedClustereasy}}
Assume for contradiction that $\beta_0$ is such that~(\ref
{eqLemmaplantedCluster2}) holds for all $k\ln2-\ln k\leq\beta\leq
\beta_0$ but $\beta_{\mathrm{crit}}(d,k)< \beta_0$. By
Corollary~\ref{Corbcrit}, we have
$\beta_{\mathrm{crit}}(d,k)\ge k\ln2-\ln k$. We pick and fix a
number\break $\beta_{\mathrm{crit}}(d,k)<\beta<
\beta_0$. We let $A=\ln2 +\frac{d}k \ln(1-2^{1-k}(1-\exp(-\beta
)))$. There exists $\eps>0$ such that
%
%
\begin{equation}
\label{eqbetastar} \lim_{n \to\infty}\frac{1}n \Erw\bigl[\ln
Z_{\beta}H_k(n,m)\bigr]< A-\eps.
\end{equation}
On the other hand, (\ref{eqLemmaplantedCluster2}) and Lemma~\ref
{LemmaMcDiarmid} ensure that we can apply Lemma~\ref
{LemmatamePlanting} and find a number $c>0$ such that
%
%
\begin{eqnarray}
\label{eqtame} &&\Erw\bigl[Z_{\beta,\mathrm{tame}}\bigl
(H_k(n,m)\bigr)\bigr]
\geq c\cdot\Erw\bigl[Z_{\beta
}\bigl(H_k(n,m)\bigr)\bigr].
\end{eqnarray}

Hence, $\Erw[Z_{\beta,\mathrm{tame}}(H_k(n,m))^2]=O(\Erw[Z_{\beta
,\mathrm{tame}}(H_k(n,m))]^2)$ by Lemma~\ref{Lemmasmm}.
Using the Paley--Zygmund inequality, there is a number $C>0$ such that
\[
\liminf_{n \to\infty}\pr\brkk{Z_{\beta,\mathrm
{tame}}\bigl(H_k(n,m)
\bigr)\geq\Erw\bigl[Z_{\beta,\mathrm{tame}}\bigl(H_k(n,m)\bigr
)\bigr]/2}
\geq1/C>0.
\]
With (\ref{eqtame}) and because $c/2\cdot\Erw[Z_{\beta}(H_k(n,m)
)]>\exp(nA-n\eps/3)$ we see that
\[
\liminf_{n \to\infty}\pr\brkk{Z_{\beta,\mathrm
{tame}}\bigl(H_k(n,m)
\bigr)\geq\exp(nA-n\eps/3)}>0.
\]
With Lemma~\ref{LemmaZAzuma}, it follows that
\[
\lim_{n\to\infty}\pr\brkk{Z_{\beta,\mathrm{tame}}\bigl(H_k(n,m)
\bigr)\geq\exp(nA-2n\eps/3)}=1.
\]
With (\ref{eqbetastar}), we get the contradiction
\[
A-\eps>\liminf_{n\to\infty}\frac{1}n \Erw\bigl[\ln
Z_{\beta,\mathrm
{tame}}\bigl(H_k(n,m)\bigr)\bigr]\geq A-2\eps/3
\]
which refutes our assumption that $\beta_{\mathrm{crit}}(d,k)< \beta_0$.
\end{pf*}

\begin{pf*}{Proof of Proposition~\ref{LemmaplantedCluster}}
The proposition is immediate from Corollary~\ref{CorollaryantiPlanting}
combined with Lemma~\ref{Lemmapartition} and from Corollary~\ref
{LemmaplantedClustereasy}.
\end{pf*}

\section{The cluster size}\label{SecclusterSize}

\textit{In this section}, \textit{we prove Proposition}~\ref{propclustersize}.
\textit{Throughout the section}, \textit{we assume that} $d/k= 2^{k-1}\ln2+O_k(1)$ \textit{and
that} $\beta\geq k\ln2-\ln k$.

In order to analyse the cluster size, we will show that there is a
large set of vertices (the ``core'') whose value cannot be changed
without creating a large number of monochromatic edges. Hence, the
contribution of these vertices to the cluster size can be controlled.
Then we analyze the contribution of the remaining vertices.

The proof strategy broadly follows the argument for estimating the
cluster size in the ``zero temperature'' case from~\cite{Lenka}.
However, the fact that we are dealing with a finite $\beta$ causes
significant complications.
More precisely, one of the key features of the ``zero temperature''
case is the existence of ``frozen variables'',
that is, vertices that take the same color in all colorings in the cluster.
Indeed, in the zero temperature case the problem of estimating the
cluster size basically reduces to estimating the number of ``frozen variables''.
By contrast, in the case of finite $\beta$, frozen variables do not
exist. In effect, we need to take a much closer look.

We let $\bolds{\sigma}:\brk  n\rightarrow\{{\pm1} \}$ be a map
chosen uniformly at
random conditioned on the event that $\bolds{\sigma}\in\Bal$ and
$\G$ be
the random hypergraph obtained by inserting each edge that is
monochromatic under $\bolds{\sigma}$ with probability $p_1$ and each
edge that
is bichromatic with probability $p_2$.

We say that a vertex $v$ \textit{supports} an edge $e \ni v$ under
$\bolds
{\sigma}
$ if $\bolds{\sigma}(e\setminus\{{v} \})=\{-\bolds{\sigma}(v)\}$.
In this case, we call $e$ \textit{critical}.
Moreover, if $U\subset\brk n$, then we say that an edge $e$ of
$\mathbf H$ is $U$-\textit{endangered} if
$\vert\bolds{\sigma}(U\cap e)\vert=1$ (i.e., the vertices in
$U\cap e$ all have the
same color).

For the first three subsections of this section, it will be convenient
to introduce a slightly more general construction. Let $\omega\geq0$
be fixed and let $v_1,\ldots,v_\omega$ be vertices chosen uniformly
at random without replacement from all vertices in $\G$.
Let $\G'$ be the hypergraph obtained from $\G$ by removing
$v_1,\ldots,v_\omega$ and edges $e$ involving one of these vertices.
Without loss of generality, we can assume that $\{v_1, \dots, v_\omega
\} = \{n-\omega+1, \dots, n\}$. The edge set of $\G'$ is thus $\brk
{n'}$, with $n' = n - \omega$.

\subsection{The core}
Let $\core(\mathbf H,\bolds{\sigma})$ be the maximal set $V'\subset
\brk n$
of vertices such that the following two conditions hold.
\begin{description}
\item[CR1] Each vertex $v\in V'$ supports at least 100 edges that
consist of vertices from $V'$ only.
\item[CR2] No vertex $v\in V'$ occurs in more than 10 edges that are
$V'$-endangered under~$\bolds{\sigma}$.
\end{description}
If $V',V''$ are sets that satisfy {\bf CR1}--{\bf CR2}, then so does
$V'\cup V''$.
Hence, the core is well-defined.

%
\begin{proposition} \label{Propcore}
A.a.s. $\vert\core(\G,\bolds{\sigma})\vert= n(1-\tilde O_k(2^{-k}))$.
\end{proposition}

To prove this proposition, we consider the following \textit{whitening
process} on the graph $\G'$ whose result $U$ is such that its
complement $\bar U=\brk {n'} \setminus U$ is a subset of $\core(\G
',\bolds{\sigma})$.
\begin{description}
\item[WH1] Let $W$ contain all vertices of $\G'$ that either support
fewer than 200 edges or that occur in more than 2 edges
that are monochromatic under $\bolds{\sigma}$.
\item[WH2] Let $U=W$ initially.
While there is a vertex $v\in[n']\setminus U$ such that:
\begin{itemize}
\item$v$ occurs in more than $5$ edges that are $\brk {n'}\setminus
U$-endangered and contain a vertex from $U$, or
\item$v$ supports fewer than $150$ edges containing vertices in $\brk
{n'}\setminus U$ only,
\end{itemize}
add $v$ to $U$.
\end{description}

Proposition~\ref{Propcore} will be a consequence of the following lemma,
by taking $\omega= 0$ and noticing that $\core(\G',\bolds{\sigma
})$ is a
superset of the set $\bar U$.

%
\begin{lemma}\label{Propcorelemmaaux}
Let $U$ be the outcome of the process {\bf WH1}--{\bf WH2} on $\G'$. Then
$\vert U\vert= n'\tilde O_k(2^{-k})$ a.a.s.
\end{lemma}
The rest of this subsection is dedicated to the proof of this lemma. We
first bound the size of the set $W$ generated by {\bf WH1}.

%
\begin{lemma}\label{LemmaWH1}
A.a.s. the set $W$ contains $n'\tilde O_k(2^{-k})$ vertices.
\end{lemma}
\begin{pf}
Our assumptions on $\beta$ and $d$ ensure that the number of\break 
monochromatic edges that any fixed vertex $v$ occurs in is binomially
distributed with mean $\tilde O_k(2^{-k})$.
Therefore, the probability that $v$ occurs in more than $2$
monochromatic edges is bounded by $\tilde O_k(2^{-2k})$.
Furthermore, the number of edges that $v$ supports is binomially
distributed with mean $k\ln2+O_k(1)$.
Hence, by the Chernoff bound the probability that $v$ supports fewer
than $200$ edges is bounded by $\tilde O_k(2^{-k})$.
Consequently,
%
%
\begin{equation}
\label{eqLemmaWH11} \Erw\bigl[ \vert W\vert\bigr] = n'\tilde
O_k\bigl(2^{-k}\bigr).
\end{equation}
Finally, either adding or removing a single edge from the hypergraph
can alter the size of $W$ by at most $k$.
Therefore, (\ref{eqLemmaWH11}) and Azuma's inequality imply that
$\vert W\vert= n'\tilde O_k(2^{-k})$ a.a.s., as desired.
\end{pf}

In the next step, we state two results excluding some properties of
small sets of vertices in $\G'$.

%
\begin{lemma}\label{LemmaWH2a}
A.a.s. the random hypergraph $\G'$ enjoys the following property:
%
%
\begin{equation}
\label{eqLemmaWH2a} \quad\parbox{10.5cm} {There is no set $T\neq\varnothing
$ of vertices
with $\vert T\vert\leq n'/k^8$ such that at least
$0.9\vert T\vert$ vertices from $T$ occur in two or more $
\brk {n'}\setminus T$-endangered edges that contain another vertex
from $T$.}
\end{equation}
\end{lemma}
\begin{pf}
For a set $T\subset\brk {n'}$ we define $\eps= \vert T\vert/n'$
and we let
$X_i(T)$ for $i\in\{2,\ldots, k\}$ be the number of edges that are
$\brk {n'}\setminus T$-endangered and contain exactly $i$ vertices from
$T$. Then $X_i(T)$ is stochastically dominated by a binomial random
variable $\Bin((1+o(1))2^{i+1-k}\bink{\eps n'}{i}\bink
{n'}{k-i},2p )$.
Indeed, there are $ \bink{\eps n'}{i}$ ways to choose $i$ vertices
from $T$ and at most $\bink{(1-\eps)n'}{k-i} \le\bink{n'}{k-i}$
ways to choose $k-i$ vertices from $\brk {n'}\setminus T$.
Moreover, these $k-i$ vertices are required to have the same color and
because we assumed that $\bolds{\sigma}$ is balanced, this gives rise
to the
$(1+o(1))2^{i+1-k}$-factor.
Let $X(T)=\sum_{i=2}^k X_i(T)$ be the total number of edges that are
$\brk {n'}\setminus T$-endangered and contain at least two vertices
from $T$. Then using the rough upper bound $\bink nk2p\leq n2^k\ln2$
we obtain
%
\begin{equation}
\label{equpperboundErwXT} \Erw\bigl[X(T)\bigr]=\sum_{i=2}^k
\Erw\bigl[X_i(T)\bigr] \le k \Erw\bigl[X_2(T)\bigr]
\le3.6 k^3 \eps^2 n'.
\end{equation}
Let ${{\mathcal E}}(T)$ be the event that $X(T)\ge1.8\vert T\vert$.
If the set $T$
satisfies (\ref{eqLemmaWH2a}) then ${{\mathcal E}}(T)$ occurs.
The Chernoff bound (Lemma~\ref{lemchernoffsum})
and the above upper bound (\ref{equpperboundErwXT}) on $\Erw
[X(T)]$ yield
\begin{eqnarray*}
\pr\bigl[{{\mathcal E}}(T)\bigr]& \le&\exp\biggl(-1.8\eps n'\ln
\biggl(\frac{1}{2ek^3\eps} \biggr) \biggr).
\end{eqnarray*}
Hence, the probability of the event ${{\mathcal E}}$ that there is a
set $T$ of
size $\vert T\vert\leq n'/k^8$ such that ${{\mathcal E}}(T)$ occurs
is bounded by
\begin{eqnarray*}
\pr\brkk{{\mathcal E}}&\leq&\sum_{T:\vert T\vert\leq n'/k^8}\pr
\brkk{{{
\mathcal E}}(T)}\leq\sum_{1/n'\leq\eps\leq1/k^8}\binkd{n'}
{\eps n'} \exp\biggl(-1.8\eps n'\ln\biggl(
\frac{1}{2ek^3\eps} \biggr) \biggr)
\\
&\leq&\sum_{1/n'\leq\eps\leq1/k^8} \biggl({\frac{2\eul n'}{\eps
n'}}
\biggr)^{\eps
n'}\exp\biggl(-1.8\eps n'\ln\biggl(
\frac{1}{2ek^3\eps} \biggr) \biggr)
\\
&\leq&\sum_{1/n'\leq\eps\leq1/k^8} \exp\bigl(\eps n'
\bigl(5+5.6\ln(k)+0.8\ln(\eps) \bigr) \bigr) =o(1),
\end{eqnarray*}
as claimed.
\end{pf}

%
\begin{lemma}\label{LemmaWH2b}
A.a.s. the random hypergraph $\G'$ enjoys the following property:
%
%
\begin{equation}
\label{eqLemmaWH2b} \parbox{10.5cm} {There is no set $T\neq\varnothing
$ of vertices of
size $\vert T\vert\leq n'/k^6$ such that at least
$0.09\vert T\vert$ vertices from $T$ support at least 20 edges that contain
another vertex from $T$.}
\end{equation}
\end{lemma}
\begin{pf}
For a set $T\subset\brk {n'}$ and a set $Q\subset\brk T$, we let
${{\mathcal E}}
(T,Q)$ be the event that each vertex $v\in Q$ supports at least 20
edges that contain another vertex from $T$.
Let $\eps=\vert T\vert/n'$.
Then for each vertex $v$ the number $X_v$ of edges that $v$ supports
and that contain another vertex from $T$ is
stochastically dominated by a binomial random variable $\Bin
((1+o(1))2^{2-k}\eps n'\bink{n'}{k-2},p_2 )$.
Indeed, there are $\eps n'-1$ ways to choose another vertex $v'\neq v$
from $T$, and at most $\bink{n'}{k-2}$
ways to choose $k-2$ further vertices to complete the edges.
Moreover, these $k-2$ vertices are required to have color $-\bolds
{\sigma}
(v)$, and because we assumed that $\bolds{\sigma}$ is balanced this
gives rise
to the $(1+o(1))2^{2-k}$-factor.
Furthermore, the random variables $X_v$ are mutually independent,
because the edges in question
are distinct as they are supported by the distinguished vertex~$v$.
Therefore, using the rough upper bound
$\bink{n}kp_2\leq n2^k\ln2$, we obtain
%
%
\begin{eqnarray}\label{eqLemmaWH2b1}
\nonumber
\pr\brkk{{{\mathcal E}}(T,Q)}&\leq&\prod_{v\in Q}
\pr\brk{X_v\geq20}
\\
&\leq& \pr\brkkk{\Bin\biggl(\bigl(1+o(1)\bigr)2^{2-k}\eps
n'\binkd{n'} {k-2},p_2 \biggr)
\geq20}^{\vert Q\vert}
\\
&\leq&\bigl(k^2\eps\bigr)^{20\vert Q\vert}.\nonumber
\end{eqnarray}

Now, let ${{\mathcal E}}(T)$ be the event that there is a set $Q\subset
\brk T$
of size $\vert Q\vert\geq0.09\vert T\vert$ such that ${{\mathcal
E}}(T,Q)$ occurs.
Then~(\ref{eqLemmaWH2b1}) implies that
\begin{eqnarray*}
\pr\brkk{{{\mathcal E}}(T)}&\leq&2^{\vert T\vert}\bigl(k^2\vert T
\vert/n'\bigr)^{1.8\vert T\vert}.
\end{eqnarray*}
Hence, the probability of the event ${{\mathcal E}}$ that there is a
set $T$ of
size $\vert T\vert\leq n'/k^6$ such that ${{\mathcal E}}(T)$ occurs
is bounded by
\begin{eqnarray*}
\pr\brk{{\mathcal E}}&\leq&\sum_{T:\vert T\vert\leq n'/k^6}\pr
\brkk{{{
\mathcal E}}(T)}\leq\sum_{1\leq t\leq n'/k^6}\binkd{n'}t2^{t}
\bigl(k^2t/n'\bigr)^{1.8t}
\\
&\leq&\sum_{1\leq t\leq n'/k^6} \biggl({\frac{2\eul n'}{t}}
\biggr)^t\bigl(k^2t/n'\bigr)^{1.8t}
\leq\sum_{1\leq t\leq n'/k^6}\brkk{2\eul\bigl(t/n'
\bigr)^{0.8} k^{3.6}}^t=o(1),
\end{eqnarray*}
as claimed.
\end{pf}

\begin{pf*}{Proof of Lemma~\ref{Propcorelemmaaux}}
By Lemmas~\ref{LemmaWH2a} and~\ref{LemmaWH2b}, we may assume that
$\G'$ enjoys the properties~(\ref{eqLemmaWH2a}) and~(\ref{eqLemmaWH2b}).
We are going to argue that $\vert U\vert\leq k\vert W\vert$
a.a.s. Indeed, assume for contradiction that $\vert U\vert>k\vert
W\vert$ and let $U'$ be the
set obtained by {\bf WH2} when precisely $(k-1)\vert W\vert$
vertices have been
added to $U$; thus, $\vert U'\vert=k\vert W\vert$.
Then by construction each vertex $v\in U'$ has one of the following properties:
\begin{longlist}[(2)]
\item[(1)] $v$ belongs to $W$,
\item[(2)]  or $v$ occurs in two or more $\brk {n'}\setminus U'$-endangered edges,
\item[(3)]  or $v$ supports at least 20 edges that contain another vertex
from $U'$.
\end{longlist}
Let $U_0\subset U'$ be the set of all $v\in U'$ that satisfy (1), let
$U_1\subset U'\setminus U_0$ be the set of all $v \in U'\setminus U_0$
that satisfy (2) and let $U_2=U'\setminus(U_0\cup U_1)$.
There are two cases to consider.
\begin{longlist}
\item[\textit{Case} 1.] $\vert U_1\vert\geq0.9\vert U'\vert$]
then (\ref{eqLemmaWH2a}) implies that $\vert U'\vert>n'/k^8$.

\item[\textit{Case} 2.] $\vert U_1\vert<0.9\vert U'\vert$]
then $\vert U_0\vert+\vert U_2\vert\geq0.1\vert U'\vert$ and
since $\vert U_0\vert=\vert W\vert$ and $\vert U'\vert=k\vert
W\vert$ we
have $\vert U_2\vert\geq0.09\vert U'\vert$ for $k$ large enough.
Thus, (\ref{eqLemmaWH2b}) entails that $\vert U'\vert>n'/k^6$.
\end{longlist}
Hence, in either case we have $k\vert W\vert=\vert U'\vert
>n'/k^8$, and thus $\vert W\vert>n'/k^9$.
But by Lemma~\ref{LemmaWH1} we have $\vert W\vert=n'\tilde O_k(2^{-k})$
a.a.s. Thus, we conclude that $\vert U\vert\leq k\vert W\vert
=n'\tilde O_k(2^{-k})$ a.a.s.
\end{pf*}

\subsection{The backbone}

We define the \textit{backbone} $\back(\mathbf H,\bolds{\sigma})$ as
the set of
all vertices $v \in[n] \setminus\core(\mathbf H,\bolds{\sigma})$
such that
the following two conditions hold.
\begin{description}
\item[BB1] $v$ supports at least one edge $e$ such that $e\setminus
\{{v} \}\subset\core(\mathbf H,\bolds{\sigma})$.
\item[BB2] $v$ does not occur in a $ \{{v} \} \cup\core
(\mathbf
H,\bolds{\sigma})$-endangered edge.
\end{description}

Given $\G'$, we simply reconstruct $\G$ (in distribution) by adding
for each $i \in\brk \omega$ each monochromatic edge involving $v_i$
with probability $p_1$, and each bichromatic edge involving $v_i$ with
probability $p_2$. We let $\cA$ be the event that:
\begin{itemize}
\item no vertex $v\in\brk {n'}$ is incident with more than one edge
containing a vertex from $\{v_1,\ldots,v_\omega\}$, and
\item there is no edge containing two vertices from $\{v_1,\ldots
,v_\omega\}$.
\end{itemize}

With the notation from the previous subsection we let $\bar U$ be the
complement of the set of vertices produced by the whitening process\vspace*{1pt}
{\bf WH1}--{\bf WH2} applied to the hypergraph $\G'$.
We note that $\vert\bar U\vert= n'(1-\tilde O_k(2^{-k}))$ a.a.s.~by Lemma~\ref
{Propcorelemmaaux}. In addition, if $\cA$ occurs, then $\bar
U\subset\core(\G, \bolds{\sigma})$. In this case, the following lemma
states the probabilities for some events concerning the vertices $v_i,
i \in[\omega]$.

%
\begin{lemma} \label{lemmanumbers} Assume that $\mathcal{A}$ holds.
Let $l \geq0$ be fixed. Then the following statements are true for all
$i \in[\omega]$:
\begin{longlist}[(2)]
\item[(1)] The probability that $v_i$ supports exactly $l$ edges is $
(1+o(1)) \frac{\lambda^l}{l! \exp(\lambda)} $ where
\[
\lambda= \frac{d}{ 2^{k-1}-1+\exp(-\beta)} = k \ln2 + \tilde O_k
\bigl(2^{-k}\bigr). %
\]
\item[(2)] The\vspace*{1pt} probability that $v_i$ occurs in exactly $l$ monochromatic
edges is $(1+o(1)) \frac{(\lambda')^l}{l! \exp(\lambda')}$ where
$\lambda'=\tilde O_k (2^{-k})$.\vspace*{1pt}
\item[(3)] The probability that there exist exactly $l$ edges blocking $v_i$
and containing at least one vertex outside $\{v_i\} \cup\bar U$ is
$(1+o(1)) \frac{(\lambda'')^l}{l! \exp(\lambda'')}$ where $\lambda
''=\tilde O_k (2^{-k})$.
\item[(4)] The probability that exactly $l$ edges are $\{v_i\} \cup\bar
U$-endangered is $(1+o(1)) \frac{(\lambda''')^l}{l! \exp(\lambda
''')}$ where $\lambda'''=\tilde O_k (2^{-k})$.
\end{longlist}
\end{lemma}

\begin{pf}
For each $i \in[\omega]$ the number of edges that $v_i$ supports is a
binomial random variable $\Bin({ n-1 \choose k-1
}(1+o(1))2^{1-k}, p_2 )$ and the number of monochromatic edges
involving $v_i$ is a binomial random variable\break $\Bin
({n-1 \choose k-1 } (1+o(1))2^{1-k}, p_1 )$. Indeed, because we
assumed that $\bolds{\sigma}$ is balanced, there are ${ n-1 \choose k-1}
(1+o(1))2^{1-k}$ edges $e$ involving $v_i$ such that $\bolds{\sigma
}(v) = -
\bolds{\sigma}(v_i)$ [resp., $\bolds{\sigma}(v) = \bolds{\sigma
}(v_i)$] for all $v \in
e\setminus\{{v_i} \} $ and each of them is added
independently at
random with probability $p_2$ (resp., $p_1$).
Hence, the Poisson approximation of the binomial distribution shows
that the probability that $v_i$ supports precisely $l$ edges is
$(1+o(1)) \frac{\lambda^l}{l! \exp(\lambda)}$ with
\[
\lambda=\pmatrix{ n-1 \cr k-1} \frac{p_2}{2^{k-1}} = \frac{d}{
2^{k-1}-1+\exp(-\beta)},
\]
which proves assertion (1).
Moreover, since $\beta= \Omega_k( k \ln2)$ and $d = \tilde
O_k(2^k)$, the probability that $v_i$ occurs in precisely $l$
monochromatic edges is $(1+o(1)) \frac{(\lambda')^l}{l! \exp(\lambda
')}$ with
\[
\lambda'=\pmatrix{ n-1 \cr k-1}\frac{p_1}{2^{k-1}} = \lambda\tilde
O_k\bigl(2^{-k}\bigr)=\tilde O_k
\bigl(2^{-k}\bigr).
\]
This implies assertion (2).

The\vspace*{1pt} probability that in an edge blocking $v_i$ at least one of the
vertices is outside $\{v_i\} \cup\bar U$ is $\tilde O_k(2^{-k})$ by
Lemma~\ref{Propcorelemmaaux}. Using (1), the number of edges
blocking $v_i$ and containing at least one vertex outside $\{v_i\} \cup
\bar U$ is stochastically dominated by a $\Bin({ n-1 \choose
k-1}\tilde O_k(4^{-k}),p_2 )$ random variable. (3) then follows
by the Poisson approximation.

If an edge $e$ is $\{v_i\} \cup\bar U$-endangered it is either
monochromatic or such that $\vert(e \setminus\{v_i \}) \cap\bar
U\vert\le
k-2$. Given $\G'$, these two events are independent and the numbers of
edges of each type are binomially distributed. The expected number of
edges of the first type is $\tilde O_k(2^{-k})$ by (2). The expected
number of edges of the second type is $\tilde O_k(2^{-k})$ by Lemma~\ref{LemmaWH1}. Thus, (4) follows again from the Poisson approximation.
\end{pf}

\subsection{The rest}

Let $\rest(\mathbf H,\bolds{\sigma}) = [n] \setminus(\core(\mathbf
H,\bolds{\sigma}
) \cup\back(\mathbf H,\bolds{\sigma}))$.

%
\begin{proposition}\label{Proprest}
A.a.s.~$ \vert\rest(\mathbf H,\bolds{\sigma}) \vert= n 2^{-k}
(1+\tilde O_k(2^{-k}))$.
\end{proposition}

\begin{pf}
$\rest(\mathbf H,\bolds{\sigma})$ contains at least all vertices
that do not
support an edge. Because the number of edges that a vertex supports is
binomially distributed with mean $k\ln2+O_k(1)$, by the Chernoff bound
we have $\vert\rest(\mathbf H,\bolds{\sigma}) \vert\ge n 2^{-k}
(1+\tilde
O_k(2^{-k}))$ a.a.s.
Now let $Y=\rest(\mathbf H,\bolds{\sigma})$ and let $\omega=\omega
(n)$ be a
slowly diverging function.
Let $\eps=\tilde O_k(2^{-k})$.
We are going to show that
%
%
\begin{equation}
\label{eqProprest} \Erw\bigl[Y(Y-1)\cdot\cdots\cdot(Y-\omega
+1)\bigr]\leq\biggl({
\frac
{(1+\eps+ o(1))n}{2^k}} \biggr)^\omega.
\end{equation}
This bound implies the assertion; indeed,
\begin{eqnarray*}
&& \pr\brkk{Y>(1+2\eps)n2^{-k}}
\\
&&\qquad  \leq\pr\brkk{Y(Y-1)\cdot\cdots\cdot(Y-\omega+1)>\bigl(\bigl
(1+2\eps-o(1)
\bigr)n2^{-k}\bigr)^\omega}
\\
&&\qquad  \leq\frac{\Erw[Y(Y-1)\cdot\cdots\cdot(Y-\omega+1)]}{((1+2\eps
-o(1))n2^{-k})^\omega} \leq\biggl({\frac{1+\eps+ o(1)}{1+2\eps
-o(1)}} \biggr)^\omega=o(1).
\end{eqnarray*}

To prove (\ref{eqProprest}), we observe that $Y(Y-1)\cdot\cdots
\cdot(Y-\omega+1)$ is just the number
of ordered $\omega$-tuples of vertices belonging to neither the core
nor the backbone---that is, belonging to $Y$.
Hence, by symmetry and the linearity of expectation,
\[
\Erw\bigl[Y(Y-1)\cdot\cdots\cdot(Y-\omega+1)\bigr]\leq n^\omega
\pr\brk{v_1,\ldots,v_\omega\in Y}.
\]

Thus, we are left to estimate $\pr\brk {v_1,\ldots,v_\omega\in Y}$.
If $\cA$ occurs, then $\bar U\subset\core(\G, \bolds{\sigma})$.
If $\bar U
\subset\core(\G, \bolds{\sigma})$ and $v_1,\ldots,v_\omega\in
Y$, then
for any $i\in\brk \omega$ one of the following must occur.
\begin{longlist}[(1)]
\item[(1)] There is no edge blocking $v_i$ that consists of vertices in
$ \{{v_i} \} \cup\bar U$ only.
\item[(2)] $v_i$ occurs in more than $10$ edges that are $ \{
{v_i} \}\cup
\bar U$-endangered.
\item[(3)] There are at least $200$ edges blocking $v_i$ but fewer than
$100$ of them consist of vertices in $ \{{v_i} \}\cup\bar
U$ only.
\item[(4)] There are at most $200$ edges blocking $v_i$ and one edge $e$
such that $v_i\in e$ and that is $ \{{v_i} \}\cup\bar U$-endangered.
\end{longlist}
Indeed, if a vertex $v_i$ is in $\rest(\G,\bolds{\sigma})$ then it violates
one of the conditions {\bf CR1} and {\bf CR2} and one of {\bf BB1} and
{\bf BB2}. Therefore, we have to consider several cases. If $v_i$
violates {\bf BB1}, then (1) is true. If it violates {\bf CR1} and {\bf
BB2}, then either (3) or~(4) is true. If $v_i$ violates {\bf CR2} and
one of {\bf BB1} and {\bf BB2}, then (2) is true.

Let $\cB_i$ be the event that one of the above is true for $i \in\brk
\omega$. By the principle of deferred decisions, we have $\pr\brk
{\cA}=1-O(\omega^2/n)$ and, therefore, we get
\begin{eqnarray*}
\pr\brk{v_1,\ldots,v_\omega\in Y} &\leq&\pr
\brk{v_1,\ldots,v_\omega\in Y\mid\cA}+o(1) \le\pr\brkkkkp{
\bigcap_{i=1}^\omega\cB_i\Big|
\cA}+o(1).
\end{eqnarray*}
Given that there is no edge containing two vertices from $v_1,\ldots
,v_\omega$, the events $\cB_1,\ldots,\cB_\omega$ are mutually independent.
Therefore, $\pr\brk {\bigcap_{i=1}^\omega\cB_i\mid\cA}=\pr\brk
{\cB
_1\mid\cA}^{\omega}$. Given\vspace*{2pt} that $\mathcal{\cA}$ occurs, by
Lemma~\ref{lemmanumbers} the probability of event (1) is asymptotically equal to
$2^{-k} + \tilde O_k(4^{-k})$ and\vspace*{1pt} the probabilities of events (2), (3)
and (4) are asymptotically equal to $\tilde O_k(4^{-k})$.
Hence, $\pr\brk {\cB_1\mid\cA}=2^{-k}+\tilde O_k(4^{-k})$ and $ \pr
\brk {v_1,\ldots,v_\omega\in Y} \leq(2^{-k}+\tilde
O_k(4^{-k})+o(1))^\omega=((1+\eps+ o(1))2^{-k})^\omega$.
\end{pf}

We define $\free(\mathbf H,\bolds{\sigma})$ as the set of all
vertices $v\in
\rest(\mathbf H,\bolds{\sigma})$ such that $v$ occurs only in edges
$e$ such
that $e \cap\core(\G, \bolds{\sigma})$ is bichromatic.

%
\begin{proposition} \label{Propfree} A.a.s.~$\vert\rest(\mathbf
H,\bolds{\sigma}
) \setminus\free(\mathbf H,\bolds{\sigma})\vert=n\tilde O_k
(4^{-k})$. In
particular, $\vert\free(\mathbf H,\bolds{\sigma})\vert=n(2^{-k}+
\tilde O_k
(4^{-k}))$.
\end{proposition}
\begin{pf} We introduce $Y= \vert\rest(\mathbf H,\bolds{\sigma
})\setminus\free
(\mathbf H,\bolds{\sigma})\vert$ and proceed just as in the proof of
Proposition~\ref
{Proprest}. To estimate $\pr\brk {v_1,\ldots,v_\omega\in Y}$ we
observe that if $\bar U \subset\core(\G, \bolds{\sigma})$ and
$v_1,\ldots,v_\omega\in Y$ then for any $i\in\brk \omega$ one of
the following
must occur.
\begin{longlist}[(2)]
\item[(1)] There is no edge blocking $v_i$ that consists of vertices in
$ \{{v_i} \} \cup\bar U$ only and $v_i$ occurs in at least
one edge
that is $ \{{v_i} \}\cup\bar U$-endangered.
\item[(2)] $v_i$ occurs in more than $10$ edges that are $ \{
{v_i} \}\cup
\bar U$-endangered.
\item[(3)] There are at least $200$ edges blocking $v_i$ but fewer than
$100$ of them consist of vertices in $ \{{v_i} \}\cup\bar
U$ only.
\item[(4)] There are at most $200$ edges blocking $v_i$ and one edge $e$
such that $v_i\in e$ and that is $ \{{v_i} \}\cup\bar U$-endangered.
\end{longlist}
Events (2), (3) and (4) are as in the proof of Proposition~\ref{Proprest}
and their probabilities are asymptotically equal to $\tilde O_k(4^{-k})$.
By Lemma~\ref{lemmanumbers}, the probability of (1) is $\tilde
O_k(4^{-k})$ and the assertion follows.
\end{pf}

{\em In the last three subsections, we calculate the cluster size
${{\mathcal C}}
_\beta(\G,\bolds{\sigma})$ up to a small error term. We proceed by first
eliminating the contribution of the vertices in the core and in a
second step the contribution of the vertices in the backbone. Finally,
we calculate the contribution of the vertices in $\rest(\G, \bolds
{\sigma}
)$. }

\subsection{Rigidity of the core}
In the following, we let $x=k^{-5}$. We first show that the cluster of
$\bolds{\sigma}$ under $\G$ mostly consists of configurations at distance
less than $2x$ from $\bolds{\sigma}$.

\begin{lemma} \label{lemmaclusterlocalv}
A.a.s.
\[
{{\mathcal C}}_\beta(\G,\bolds{\sigma}) \sim\sum_{\tau\in\{
{-1,1} \}^n:\scal
{\bolds{\sigma}}{\tau}\geq(1-x)n}
\exp\bigl(-\beta E_{\G}(\tau)\bigr).
\]
\end{lemma}

To prove this result, we recall the notation from Section~\ref
{Secbcrit}. We need the following technical lemma.

\begin{lemma} \label{LemmaLambdabeta} Let $d/k=2^{k-1}\ln2+O_k(1)$
and $\beta\ge k\ln2-\ln k$. Then $\sup_{\alpha\in[2/3,1-k^{-5}]}
\Lambda_\beta(\alpha)<\Lambda_\beta(1)-\Omega_k(k^{-5})$.
\end{lemma}

\begin{pf}
We observe that for $\alpha\in[1-k^{-5}, 1-k^{-7}]$,
%
%
\begin{equation}
\label{eqLambdaa} \Lambda_\beta'(\alpha) =
\frac{ \ln(1- \alpha)}{2} + \frac{
d}{2^{k} } + \tilde O_k
\bigl(2^{-k}\bigr) = k \ln2 + O_k(\ln k)\geq1.
\end{equation}
An expansion of $\Lambda_\beta(\alpha)$ near $\alpha= 1$ gives
$\Lambda_\beta(1-k^{-7}) \leq\Lambda_\beta(1) + O_k(k^{-6}) $ and
together with (\ref{eqLambdaa})
this implies
%
%
\begin{equation}
\label{equpper} \Lambda_\beta\bigl(1-k^{-5}\bigr) \leq
\Lambda_\beta(1) - \Omega_k\bigl(k^{-5}\bigr).
\end{equation}
Further, using that $\Lambda'_\beta(\alpha)>0$ if $\alpha>1-1.99\ln
k /k$ (as in the proof of Lemma~\ref{Lemmasmm1}) and (\ref{equpper})
we obtain
%
%
\begin{equation}
\label{eqpart1} \sup_{\alpha\in[1-1.99 \ln k /k, 1-k^{-5}]} \Lambda
_\beta(\alpha)
\leq\Lambda_\beta\bigl(1-k^{-5}\bigr) \leq
\Lambda_\beta(1) -\Omega_k\bigl(k^{-5}\bigr).
\end{equation}
A study of $\Lambda_\beta(\alpha)$ also gives
%
%
\begin{equation}
\label{eqpart2} \sup_{\gamma\in[1.99,2.01]} \Lambda_\beta(1-
\gamma
\ln k / k) \leq\Lambda_\beta(1) - \Omega_k
\bigl(k^{-5}\bigr)
\end{equation}
and $ \Lambda_\beta(\alpha) - \Lambda_\beta(1-2.01 \ln k / k) =
\mathcal{H} ( \frac{1+\alpha}{2} ) + \tilde O_k (
(\frac{2}{2.01} )^k ) \leq0 $
for $\alpha\in[2/3, 1-2.01 \ln k /k]$, which leads to
%
%
\begin{eqnarray}
\label{eqpart3} &&\sup_{\alpha\in[2/3, 1-2.01 \ln k /k]} \Lambda
_\beta(\alpha)\nonumber
\\
&&\qquad  \leq\mathcal{H} \biggl( \frac{1+\alpha}{2} \biggr) + \tilde O_k
\biggl( \biggl(\frac{2}{2.01} \biggr)^k \biggr) +
\Lambda_\beta(1-2.01 \ln k/k)
\\
&&\qquad \leq\Lambda_\beta(1) -\Omega_k\bigl(k^{-5}
\bigr).
\nonumber
\end{eqnarray}
Combining (\ref{eqpart1}), (\ref{eqpart2}) and (\ref{eqpart3})
completes the proof of the assertion.
\end{pf}

\begin{pf*}{Proof of Lemma~\ref{lemmaclusterlocalv}}
Let $\mathcal A$ be the event that $\vert e(H_k(n,p))-m\vert\le
m^{2/3}$. Given
$\bolds{\sigma}$ and $\alpha\in[-1,1]$ and using Lemma~\ref
{Lemrelatepartition} we have
\begin{eqnarray*}
&&\Erw\brkkkkp{\sum_{\tau\in\{{-1,1} \}^n:\scal
{\bolds{\sigma}}{\tau
}= \alpha n } \exp\bigl( - \beta
E_{\G}(\tau) \bigr)\Big| \bigl\vert e(\G)-m\bigr\vert\le
m^{2/3}}
\\
&&\qquad  = \frac{ \Erw\brk{\sum_{\tau:\scal{\bolds{\sigma}}{\tau}=
\alpha n} \exp(- \beta E_{H_k(n,p)}(\bolds{\sigma}) ) \exp
(- \beta E_{H_k(n,p)}(\tau) ) \mid\mathcal A} }{ \Erw\brk
{\exp(- \beta E_{H_k(n,p)}(\bolds{\sigma}) ) \mid\mathcal A} }
\\
&&\qquad  \le\frac{\Erw[ Z_\beta(\alpha) ]}{\Erw[
Z_\beta(H'_k(n,m)) ]}\exp\bigl(O\bigl(m^{2/3}\bigr)\bigr).
\end{eqnarray*}
In order to derive the last line, we used an observation similar to
equation (\ref{eqlowerpart}) and Lemma~\ref{Lemrelatepartition}.
We observe that by Lemma~\ref{LemavgEngergy} we have a.a.s.
${{\mathcal C}}
_\beta(\G,\bolds{\sigma}) \geq\exp(- \beta E_{\G}(\bolds{\sigma
})) \sim\exp(-n
\tilde O_k(2^{-k}))$. Hence,
\begin{eqnarray*}
&&\Erw\brkkkkp {\mathop{\sum_{\tau\in\{{-1,1} \}
^n:}}_{{2/3n \leq
\scal{\bolds{\sigma}}{\tau}< (1-x)n}}\exp\bigl(-\beta
E_{\G}(\tau)\bigr)\mid\bigl\vert e(\G)-m\bigr\vert\le
m^{2/3}}
\\
&&\qquad  \le\sum_{\nu=0}^n \frac{\Erw\brk{Z_\beta(2 \nu/n-1 )}}{\Erw
\brk{Z_\beta(H'_k(n,m))}}
\mathbf{1}_{2\nu/n-1 \in[2/3,(1-x)] } \exp\bigl(O\bigl
(m^{2/3}\bigr)\bigr)
\\
&&\qquad  \leq\exp\Bigl( n \Bigl( \sup_{\alpha\in[2/3,1-x]} \Lambda
_\beta(\alpha) - \Lambda_\beta(1)+\tilde O_k
\bigl(2^{-k}\bigr) \Bigr) \Bigr){{\mathcal C}}_\beta(\G,
\bolds{\sigma})
\\
&&\qquad \le\exp\bigl(-n \Omega_k\bigl(k^{-5}\bigr)\bigr) {{
\mathcal C}}_\beta(\G,\bolds{\sigma})
\end{eqnarray*}
by Lemma~\ref{LemmaLambda} and by Lemma~\ref{LemmaLambdabeta}.
It follows from Markov's inequality that a.a.s.
\[
\sum_{\tau\in\{{-1,1} \}^n:2/3n \leq\scal{\bolds{\sigma}
}{\tau}<
(1-x)n}\exp\bigl(-\beta E_{\G}(\tau)
\bigr)=o\bigl( {{\mathcal C}}_\beta(\G,\bolds{\sigma})\bigr).
\]\upqed
\end{pf*}

We now approximate ${{\mathcal C}}_\beta(\G,\bolds{\sigma})$ based
on the previous
decomposition of the vertex set $V$. Given a $k$-uniform hypergraph $\G
$, $\bolds{\sigma}: [n] \to\{\pm1\}$, and three maps $\tau_\core:
\core
(\G,\bolds{\sigma}) \to\{ \pm1\}$, $\tau_\back: \back(\G
,\bolds{\sigma}) \to
\{ \pm1\}$ and $\tau_\rest: \rest(\G,\bolds{\sigma}) \to\{ \pm
1\}$, we
define $E_{\G}(\tau_\core, \tau_\back, \tau_\rest)$ as $E_{\G
}(\tau)$ for the unique $\tau$ whose restriction to $\core(\G
,\bolds{\sigma})$ [resp., $\back(\G,\bolds{\sigma}), \rest(\G
,\bolds{\sigma})$] is given
by $\tau_\core$ (resp., $\tau_\back, \tau_\rest$).

We introduce the ``restricted'' cluster size
\[
{{\mathcal C}}_\beta^{\back+\rest}(\G, \bolds{\sigma}) = \sum
_{\tau
_\back, \tau
_\rest} \exp\bigl( - \beta E_{\G}(
\bolds{\sigma}_\core,\tau_\back,\tau_\rest) \bigr).
\]
The summation is over $\tau_\back: \back(\G,\bolds{\sigma}) \to
\{ \pm1\}
$ and $\tau_\rest: \rest(\G,\bolds{\sigma}) \to\{ \pm1\}$.
The aim of this section is to prove the following.

\begin{proposition} \label{Lemmacorerigid}
A.a.s.
\[
\frac{1}{n} \ln{{\mathcal C}}_\beta^{\back+\rest}(\G,\bolds
{\sigma})
\leq\frac
{1}{n} \ln{{\mathcal C}}_\beta(\G,\bolds{\sigma}) \leq
\frac{1}{n} \ln{{\mathcal C}}_\beta^{\back+\rest}(\G,\bolds
{\sigma}) +
\exp(- 88\beta).
\]
\end{proposition}

In order to proceed, we first need a few additional results. We
introduce the set ${{\mathcal E}}_{\G}(\tau,\bolds{\sigma})$ of
edges that:
\begin{itemize}
\item are supported by a vertex $v$ such that $\tau_\core(v) \neq
\bolds{\sigma}_\core(v)$,
\item contain two or more vertices $v'$ such that $\tau_\core(v')
\neq\bolds{\sigma}_\core(v')$.
\end{itemize}
The following lemma is reminiscent of \cite{Lenka}, Lemma~5.9.

%
\begin{lemma} \label{Lemmacorecritical}
A.a.s.~it holds that, for all $\tau: [n] \to\{ \pm1\}$ satisfying
$\langle\bolds{\sigma}, \tau\rangle\geq(1-x)n$,
\[
\bigl\vert{{\mathcal E}}_{\G}(\tau,\bolds{\sigma})\bigr\vert
\leq2 \bigl
\vert\bigl\{ v: \bolds{\sigma}_\core(v) \neq\tau_\core(v)\bigr
\}
\bigr\vert.
\]
\end{lemma}

\begin{pf}
We claim that a.a.s.~$\G$ has the following property. Let $T \subset
V$ be of size $ \vert T\vert\leq n/(2 e^3 k^2 \lambda^2)$. Then
there are no
more than $2\vert T\vert$ edges that are supported by a vertex in
$T$ and
contain a second vertex from $T$. Indeed, by a first moment argument,
with $\vert T\vert=tn$ the probability that there is a set $T$ that
violates the
above property is bounded by
\begin{eqnarray*}
\pmatrix{n \cr tn} \pmatrix{ \bigl(1+o(1)\bigr) \lambda n \cr 2 t n}
\bigl(kt^2\bigr)^{2tn}& \leq&\biggl[ \bigl(1+o(1)\bigr)
\frac{e}{t} \biggl(\frac{\lambda e }{2 t} \biggr)^2 \bigl( k
t^2 \bigr)^2 \biggr]^{tn}
\\
& \leq&\bigl( \bigl(1+o(1)\bigr) t \bigl( e^3 \lambda^2
k^2 \bigr) \bigr)^{tn} = o(1).
\end{eqnarray*}

With $T = \{v: \bolds{\sigma}_\core(v) \neq\tau_\core(v) \}$ and
$x=k^{-5}$, we have $ \vert T\vert\leq2 xn <n/(2 e^3 k^2 \lambda
^2)$ which
completes the proof.
\end{pf}

%
\begin{lemma} \label{Lemmacoreenergy}
A.a.s. it holds that, for all $\tau: [n] \to\{ \pm1 \}$ satisfying
$\langle\bolds{\sigma}, \tau\rangle\geq(1-x)n$,
\[
E_{\G}(\tau_\core,\tau_\back,
\tau_\rest) \geq E_{\G}(\bolds{\sigma}_\core,
\tau_\back,\tau_\rest) + 88 \operatorname{dist}(\tau_\core,
\bolds{\sigma}_\core). %
\]
\end{lemma}

\begin{pf}
Denote for a vertex $v \in V$ and $\tau: [n] \to\{ \pm1 \}$ by:
\begin{itemize}
\item$X(v)$ the number of {critical} (under $\bolds{\sigma}$) edges $e$
supported by $v$ such that $e \setminus\{v\} \subset\core(\G,
\bolds{\sigma})$,
\item$Y(v)$ the number of $\core(\G,\bolds{\sigma})$-endangered edges
containing $v$,
\item$M_\tau(v)$ the number of edges containing $v$ that are
monochromatic under $(\bolds{\sigma}_\core,\tau_\back, \tau_\rest)$.
\end{itemize}
We can lower bound $E_{\G}(\tau_\core,\tau_\back, \tau_\rest)$
in terms of $E_{\G}(\bolds{\sigma}_\core,\tau_\back, \tau_\rest
)$ as
\begin{eqnarray}\label{eqauxcorerigidity2}
\qquad E_{\G}(\tau_\core,\tau_\back,
\tau_\rest)
&\geq& E_{\G}(\bolds{\sigma}_\core,\tau_\back,
\tau_\rest)
\nonumber\\[-8pt]\\[-8pt]\nonumber
&&{} + \sum_{v:\tau_\core(v) \neq\bolds{\sigma}_\core
(v)} \bigl(X(v) -
M_\tau(v)\bigr) - \bigl\vert{{\mathcal E}}_{\G}(\tau,
\bolds{\sigma})\bigr\vert.
\nonumber
\end{eqnarray}
Only edges that were $\core(\G,\bolds{\sigma})$-endangered can be
monochromatic under $(\bolds{\sigma}_\core,\tau_\back,\tau_\rest
)$: $
M_{\tau}(v) \leq Y(v)$. In particular,
%
\begin{equation}
\label{eqauxcorerigidity1} \forall v \in\core(\G,\bolds{\sigma}),\qquad X(v) -
M_\tau(v)
\geq90.
\end{equation}
On the other hand, we can upper bound $\vert{{\mathcal E}}_{\G}(\tau
,\bolds{\sigma})\vert$
with Lemma~\ref{Lemmacorecritical}. Replacing in (\ref
{eqauxcorerigidity2}) and using (\ref{eqauxcorerigidity1}) gives
\[
E_{\G}(\tau_\core,\tau_\back,
\tau_\rest) \geq E_{\G}(\bolds{\sigma}_\core,
\tau_\back, \tau_\rest) + 88\operatorname{dist}(
\tau_\core, \bolds{\sigma}_\core), %
\]
a.a.s., completing the proof.
\end{pf}

\begin{pf*}{Proof of Proposition~\ref{Lemmacorerigid}}
We first prove the lower bound on ${{\mathcal C}}_\beta(\G,\bolds
{\sigma}
)$. With Proposition~\ref{Propcore}, a.a.s.~for all $(\tau_\back,
\tau
_\rest)$ we have $\langle\bolds{\sigma}, (\bolds{\sigma}_\core
,\tau
_\back,\tau_\rest) \rangle\geq(1-x)n$. Hence, with Lemma~\ref
{lemmaclusterlocalv}. a.a.s.
\[
{{\mathcal C}}_\beta(\G,\bolds{\sigma}) \geq\sum_{\tau_\back,
\tau
_\rest}
\exp\bigl(- \beta E_{\G}(\bolds{\sigma}_\core,\tau_\back,
\tau_\rest)\bigr) = {{\mathcal C}} _\beta^{\back+\rest}(\G,
\bolds{\sigma}).
\]
To derive the upper bound, we write
%
%
\begin{eqnarray}\label{eqextraeq64v1}
{{\mathcal C}}_\beta(\G,\bolds{\sigma}) &\leq&\mathop{\sum
_{ \tau_\core:}}_{
\langle\bolds{\sigma}_\core, \tau_\core\rangle\geq(1-x)n } \sum
_{\tau_\back, \tau_\rest}
\exp\bigl(- \beta E_{\G}(\tau_\core,\tau_\back,
\tau_\rest)\bigr)
\nonumber\\[-8pt]\\[-8pt]\nonumber
&\leq& \mathop{\sum_{\tau_\core:}}_{\langle
\bolds{\sigma}_\core, \tau_\core\rangle\geq(1-x)n }
\exp\bigl(- 88\beta\operatorname{dist}(\bolds{\sigma}_\core,\tau_\core
)\bigr)
{{\mathcal C}}_\beta^{\back
+\rest}(\G,\bolds{\sigma}),
\end{eqnarray}
where the second inequality holds a.a.s.~by Lemma~\ref
{Lemmacoreenergy}. Finally,
%
%
\begin{eqnarray}
&&\mathop{\sum_{\tau_\core:}}_{\langle\bolds{\sigma}_\core, \tau
_\core
\rangle\geq(1-x)n } \exp\bigl(- 88
\beta\operatorname{dist}(\bolds{\sigma}_\core\tau_\core)\bigr)
\nonumber
\\
&&\qquad = \sum_{i=0}^{xn /2} \pmatrix{n \cr i} \exp(- 88
\beta i) \leq\sum_{i=0}^{n} \pmatrix{n \cr i}
\exp(-88\beta i)
\\
&&\qquad = \bigl(1+\exp(-88\beta)\bigr)^n \leq\exp\bigl(n \exp(-
88\beta)
\bigr).\nonumber \label{eqextraeq64v2}
\end{eqnarray}
Replacing with (\ref{eqextraeq64v2}) in (\ref
{eqextraeq64v1}) completes the proof.
\end{pf*}

\subsection{Rigidity of the backbone}
We proceed one step further by eliminating the vertices in the backbone
and comparing ${{\mathcal C}}_\beta^{\back+\rest}(\G,\bolds{\sigma
})$ to
${{\mathcal C}}_\beta
^{\rest}(\G,\bolds{\sigma})$, where
\[
{{\mathcal C}}_\beta^{\rest}(\G,\bolds{\sigma}) = \sum
_{\tau_\rest} \exp\bigl(- \beta E_{\G}(
\bolds{\sigma}_\core, \bolds{\sigma}_\back, \tau_\rest)\bigr).
\]
The sum is over $\tau_\rest: \rest(\G, \bolds{\sigma}) \to\{ \pm
1\}$. We
prove the following result.

\begin{proposition} \label{Proprigiditybackbone} A.a.s.
\[
\frac{1}{n} \ln{{\mathcal C}}_\beta^{\rest}(\G,\bolds{\sigma})
\leq\frac{1}{n} \ln{{\mathcal C}}_\beta^{\back+\rest}(\G,\bolds
{\sigma})
\leq\frac
{1}{n} \ln{{\mathcal C}}_\beta^{\rest}(\G,
\bolds{\sigma})+ \tilde O_k\bigl(4^{-k}\bigr).
\]
\end{proposition}
\begin{pf} The left inequality is obvious. To prove the right
inequality, we observe that, by definition of the backbone, for any
$\tau_\back: \back(\G, \bolds{\sigma}) \to\{ \pm1\}$ and $\tau
_\rest:
\rest(\G, \bolds{\sigma}) \to\{ \pm1\}$, the following is true.
%
%
\begin{equation}
\label{eqauxrigiditybackbone}
\qquad E_{\G}(\bolds{\sigma}_\core, \tau
_\back,
\tau_\rest) \geq E_{\G}(\bolds{\sigma}_\core,
\bolds{\sigma}_\back, \tau_\rest) + \operatorname{dist}(\bolds{\sigma
}_\back,
\tau_\back).
\end{equation}
Indeed for any vertex $v \in\back(\G, \bolds{\sigma}) $ with
$\bolds{\sigma}_\back
(v)\ne\tau_\back(v)$ and any edge $e \ni v$:
\begin{itemize}
\item either $v$ supports $e$ and $e \setminus\{v\} \subset\core(\G
,\bolds{\sigma})$, in which case the edge $e$ is bichromatic under
$(\bolds{\sigma}
_\core,\bolds{\sigma}_\back,\break \tau_\rest)$ and monochromatic under
$(\bolds{\sigma}
_\core,\tau_\back,\tau_\rest)$,
\item or $e$ is not $\{v\} \cup\core(\G,\bolds{\sigma
})$-endangered and is
bichromatic both under $(\bolds{\sigma}_\core, \bolds{\sigma}_\back
,\tau_\rest)$
and under $(\bolds{\sigma}_\core,\tau_\back,\tau_\rest)$.
\end{itemize}
Moreover, by the definition of $\back(\G, \bolds{\sigma})$ there is
at least
one edge of the first type for any $v \in\back(\G, \bolds{\sigma})
$ with
$\bolds{\sigma}_\back(v)\ne\tau_\back(v)$.\vspace*{1pt}

Using the definition of ${{\mathcal C}}_\beta^{\back+\rest}(\G
,\bolds{\sigma})$ and
(\ref{eqauxrigiditybackbone}) yields
%
%
\begin{eqnarray}
\label{l1}
&& {{\mathcal C}}_\beta^{\back+\rest}(\G,\bolds{\sigma})\nonumber
\\
&&\qquad \leq\sum_{\tau_\back, \tau_\rest} \exp\bigl(- \beta\operatorname{dist}(
\bolds{\sigma}_\back, \tau_\back)\bigr) \exp\bigl(- \beta
E_{\G}(\bolds{\sigma}_\core,\bolds{\sigma}_\back,
\tau_\rest)\bigr)
\\
&&\qquad \leq\sum_{\tau_\back} \exp\bigl(- \beta\operatorname{dist}(
\bolds{\sigma}_\back, \tau_\back)\bigr) {{\mathcal
C}}_\beta^{\rest}(\G,\bolds{\sigma}).
\nonumber
\end{eqnarray}
The remaining sum can easily be upper-bounded:
%
%
\begin{eqnarray}
\label{l2} \sum_{\tau_\back} \exp\bigl(- \beta
\operatorname{dist}(\bolds{\sigma}_\back, \tau_\back)\bigr)&=& \sum
_{i=0}^{\vert\back(\G,\bolds{\sigma})\vert} \pmatrix{\bigl\vert\back
(\G,\bolds{\sigma})\bigr
\vert\cr i} \exp(-\beta i)
\nonumber
\\
&=& \bigl(1+\exp(-\beta)\bigr)^{\vert\back(\G,\bolds{\sigma
})\vert}
\\
& \leq& \exp\bigl(\exp(- \beta) \bigl\vert\back(\G,\bolds
{\sigma})\bigr\vert\bigr).\nonumber
\end{eqnarray}
The upper bound of Proposition~\ref{Proprigiditybackbone} then follows
from (\ref{l1}) and (\ref{l2}) combined with Proposition~\ref{Propcore}.
\end{pf}

\subsection{The remaining vertices}
We finally deal with the vertices that belong neither to the core nor
to the backbone. As anticipated in Proposition~\ref{Propfree}, most of
them are free. This yields the following result.

%
\begin{proposition} \label{Lemmarigidityrest} A.a.s.
\[
\frac{1}{n} \ln{{\mathcal C}}_\beta^{\rest}(\G,\bolds{\sigma}) =
\frac
{\ln2}{2^k} - \beta\frac{E_{\G}(\bolds{\sigma})}{n}+ \tilde O_k
\bigl(4^{-k}\bigr).
\]
\end{proposition}

In order to prove the proposition, we need the following result. Let
$M'_{\bolds{\sigma}}(v)$ be the number of monochromatic edges
involving $v$ in
the configuration $\bolds{\sigma}$.

%
\begin{lemma} \label{Lemmaenergystronglip} A.a.s.
\[
\sum_{v \in\rest(\G,\bolds{\sigma}) \setminus\operatorname{free}(\G
,\bolds{\sigma})} M'_{\bolds{\sigma}}(v) = n
\tilde O_k\bigl(4^{-k}\bigr).
\]
\end{lemma}

\begin{pf}
We start with the following observation:
\begin{eqnarray*}
&&\sum_{v \in\rest(\G,\bolds{\sigma}) \setminus\operatorname{free}(\G
,\bolds{\sigma})} M'_{\bolds{\sigma}}(v) \leq
\sum_{v: M'_{\bolds{\sigma}}(v) >2} M'_{\bolds{\sigma}}(v) + 2
\bigl\vert\rest(\G,\bolds{\sigma}) \setminus\operatorname{free}(\G
,\bolds{\sigma})\bigr\vert.
\end{eqnarray*}
The number of monochromatic edges involving a vertex $v$ is a binomial
random variable $\Bin({n-1 \choose k-1 } (1+o(1))2^{k-1},
p_1 )$. Hence $\sum_{v \in V: M'_{\bolds{\sigma}}(v) >2} M'_{\bolds
{\sigma}
}(v) = n\tilde O_k(4^{-k})$. Applying Proposition~\ref{Propfree} completes
the proof.
\end{pf}

\begin{pf*}{Proof of Proposition~\ref{Lemmarigidityrest}}
By the definition of ${\free}(\G,\bolds{\sigma})$, the number of
monochromatic edges $E_{\G}(\bolds{\sigma}_\core, \bolds{\sigma
}_\back,\tau_\rest
)$ does not depend on the values $\tau_\rest(v)$ for $v \in\free(\G
,\bolds{\sigma})$. Consequently,
\[
{{\mathcal C}}_\beta^{\rest}(\G,\bolds{\sigma}) \geq2^{\vert
\operatorname{free}(\G
,\bolds{\sigma})\vert}
\exp\bigl(-\beta E_{\G}(\bolds{\sigma})\bigr).
\]
Together with Proposition~\ref{Propfree} this gives the lower bound on
$\frac{1}{n} \ln{{\mathcal C}}_\beta^{\rest}(\G,\bolds{\sigma})$.
For the upper bound, we start with the general inequality
\[
\frac{1}{n} \ln{{\mathcal C}}_\beta^{\rest}(\G,\bolds{\sigma})
\leq\frac
{\ln
2}{n} \bigl\vert\rest(\G,\bolds{\sigma})\bigr\vert-
\frac{\beta}{n} \inf_{\tau_\rest} E_{\G}(
\bolds{\sigma}_\core, \bolds{\sigma}_\back, \tau_\rest).
\]
Because the number of monochromatic edges does not depend on the values
of the vertices in $\operatorname{free}(\G,\bolds{\sigma})$ we have
\[
\inf_{\tau_\rest} E_{\G}(\bolds{\sigma}_\core,
\bolds{\sigma}_\back, \tau_\rest) \geq E_{\G}(\bolds{\sigma}) -
\sum_{v \in\rest(\G,\bolds{\sigma})
\setminus\free(\G, \bolds{\sigma})} M'_{\bolds{\sigma}}(v).
\]
Hence, we obtain
%
%
\begin{eqnarray}\label{eqauxremaining}
&& \frac{1}{n} \ln{{\mathcal C}}_\beta^{\rest}(\G,\bolds{\sigma})
\nonumber\\[-8pt]\\[-8pt]\nonumber
&&\qquad  \leq\frac{\ln2}{n} \bigl\vert\rest(\G,\bolds{\sigma})\bigr
\vert- \beta
\frac{ E_{\G
}(\bolds{\sigma})}{n}+ \frac{\beta}{n} \sum_{v \in\rest(\G
,\bolds{\sigma})
\setminus\operatorname{free}(\G,\bolds{\sigma})}
M'_{\bolds{\sigma}}(v).
\end{eqnarray}
The upper bound follows by combining (\ref{eqauxremaining}) with
Proposition~\ref{Proprest} and Lem\-ma~\ref{Lemmaenergystronglip}.
\end{pf*}

\subsection{Proof of Proposition~\texorpdfstring{\protect\ref{propclustersize}}{3.5}}
Combining Propositions~\ref{Lemmacorerigid},~\ref
{Proprigiditybackbone} and \ref{Lemmarigidityrest}, we obtain that a.a.s.
%
%
\begin{equation}
\label{eqproofboundsclustersend} \frac{1}{n} \ln{{\mathcal C}}_\beta
(\G,\bolds{\sigma}) =
\frac{\ln
2}{2^k}- \beta\frac{E_{\G}(\bolds{\sigma})}{n} + \tilde O_k
\bigl(4^{-k}\bigr).
\end{equation}
The number of monochromatic edges in the planted model is tightly
concentrated by Chernoff bounds. Therefore, we get a.a.s.
\begin{eqnarray*}
&&E_{\G}(\bolds{\sigma}) = \pmatrix{n
\cr
k}2^{1-k}p_1
\bigl(1+o(1)\bigr) \sim\frac{\exp(-\beta)}{2^{k-1}-1+\exp(-\beta
))}\frac{d}{k}n.
\end{eqnarray*}
For $ d/k = 2^{k-1}\ln2 + O_k(1)$ and $\beta\geq k \ln2 - \ln k$, we
have $E_{\G}(\bolds{\sigma})= \break \ln2\exp(-\beta)n + \tilde
O_k(4^{-k}) n$. Inserting this in (\ref{eqproofboundsclustersend})
yields a.a.s.
\[
\frac{1}{n} \ln{{\mathcal C}}_\beta(\G,\bolds{\sigma}) =
\frac{\ln
2}{2^k} - \beta\ln2 \exp(- \beta) + \tilde O_k
\bigl(4^{-k}\bigr),
\]
proving Proposition~\ref{propclustersize}.




%

\printaddresses
\end{document}